\documentclass[11pt]{article}
\usepackage{amsmath}
\usepackage{amssymb}
\usepackage{amscd}
\usepackage{amsthm}

\newtheorem{theorem}{Theorem}[section]
\newtheorem{main}{Theorem}
\newtheorem{proposition}[theorem]{Proposition}
\newtheorem{lemma}[theorem]{Lemma}

\newtheorem{corollary}[theorem]{Corollary}
\newtheorem{conjecture}[theorem]{Conjecture}
\newtheorem{D}[theorem]{Definition}
\newenvironment{definition}{\begin{D} \rm }{\end{D}}
\newtheorem{R}[theorem]{Remark}
\newenvironment{remark}{\begin{R}\rm }{\end{R}}

\def\Cee{\mathbb{C}}
\def\Pee{\mathbb{P}}
\def\Id{\operatorname{Id}}
\def\Pic{\operatorname{Pic}}
\def\Spec{\operatorname{Spec}}
\def\Ext{\operatorname{Ext}}
\def\Ker{\operatorname{Ker}}
\def\scrO{\mathcal{O}}

\title{Cubic Threefolds and Abelian Varieties of Dimension Five}
\author{Sebastian Casalaina-Martin\thanks{The first author was partially
    supported by a VIGRE fellowship from NSF grant DMS-98-10750.}
\  and Robert Friedman\thanks{The second author was partially supported
   by NSF grant DMS-02-00810.}}

\begin{document}

\maketitle

\section*{Introduction}

Cubic threefolds have been studied in algebraic geometry since classical times. In 
\cite{CG}, Clemens and Griffiths proved that the intermediate Jacobian $JX$ of a
smooth cubic threefold is not isomorphic as a principally polarized abelian variety to
a product of Jacobians of curves,  which implies the irrationality of $X$. They also
established the Torelli theorem for cubic threefolds: the  principally polarized abelian variety
$JX$ determines $X$ up to isomorphism. A more precise result was stated by Mumford in
\cite{Mumf2}:

\begin{main} Let $X$ be a smooth cubic threefold in $\Pee^4$ with intermediate Jacobian $JX$.
Then the theta divisor of the $5$-dimensional principally polarized abelian variety $JX$ has a
unique singular point, which has multiplicity three. Moreover, the tangent cone to the theta
divisor at the singular point is a cubic hypersurface in $\Pee^4$ which is projectively
equivalent to $X$
\end{main}

Detailed proofs of Mumford's theorem have been given by R. Smith (see for example
\cite[Example 4.3.2, p.\ 80]{Fulton}) and A. Beauville \cite{Beau2.5}.

Our goal in this paper is to prove a converse to Mumford's theorem:

\begin{main}\label{ourthm} Let $A$ be a principally polarized abelian variety of dimension $5$.
Suppose that the theta divisor $\Theta$ of $A$ has a unique singular point $p$, with
$\operatorname{mult}_p\Theta = 3$. Then there exists a smooth cubic threefold $X$ such that
$A\cong JX$ as principally polarized abelian varieties.
\end{main}

To prove Theorem~\ref{ourthm}, we use the well-known connection between $JX$ and Prym
varieties, first established by Mumford \cite{Mumf2}. Let $X$ be a smooth cubic threefold, and
let $\lambda$ be a line contained in $X$. Projecting $X$ from $\lambda$ into $\Pee^2$ exhibits
the blowup of $X$ along $\lambda$ as a conic bundle over $\Pee^2$, whose discriminant curve is
a plane quintic $C$, smooth for a generic choice of $\lambda$, with an odd theta
characteristic, corresponding to a double cover $\pi\colon \widetilde{C} \to C$. Moreover,
$JX$ is isomorphic to the corresponding Prym variety $P(\widetilde{C}, \pi)$. From
this, it is easy to determine the singularities of the theta divisor of $JX$.

To prove the converse, one has to be able to analyze the singularities of the theta divisor of
an arbitrary principally polarized abelian variety of dimension $5$. The general such abelian
variety is the Prym variety of a smooth curve of genus $6$, and every principally polarized
abelian variety of dimension $5$ is the generalized Prym variety $P(\widetilde{C},
\pi)$ corresponding to an admissible double cover of a nodal curve in the sense of
Beauville \cite{Beau1}.

Thus the problem is connected to the problem of describing the singularities of the theta
divisor of a Prym variety $P(\widetilde{C}, \pi)$. The set-theoretic description of
these singularities goes back to Mumford \cite{Mumf2}. The Prym variety
$P(\widetilde{C}, \pi)=P$ may be viewed as a subvariety of the Jacobian
$J\widetilde{C}$, and the intersection of the theta divisor
$\widetilde{\Theta}$ of $J\widetilde{C}$ with $P$ is equal to $2\Xi$, where $\Xi$ is a
principal  polarization on $P(\widetilde{C}, \pi)$. From this it follows that a point $x\in P$
lies in the singular locus of
$\Xi$ if and only if either $\operatorname{mult}_x\widetilde{\Theta} \geq 4$, or
$\operatorname{mult}_x\widetilde{\Theta} = 2$ but the tangent space to $P$ at $x$ is
contained in the tangent cone to $\widetilde{\Theta}$ at $x$. These 
singularities are called \textsl{stable} and \textsl{exceptional} respectively. However, as
pointed out by Smith and Varley \cite{SV}, it is better to distinguish instead the two
cases: (1) the  tangent space to $P$ at $x$ is not contained in the tangent cone to
$\widetilde{\Theta}$ at $x$ and (2) the tangent space to $P$ at $x$ is
contained in the tangent cone to $\widetilde{\Theta}$ at $x$. In the first case, the
multiplicity of $\Xi$ at $x$ is exactly one half the multiplicity of $\widetilde{\Theta}$ at
$x$, whereas in the second case there is only an inequality. One goal of this paper is to
develop tools to handle the case (2). While we only discuss here some of the
applications to cubic threefolds, the techniques of this paper can be developed much
further. For example, one can prove directly that the multiplicity of the
singular point of the theta divisor of $JX$ is three and use this to give another
proof of the Torelli theorem for cubic threefolds. It is possible to analyze the
singularities of Prym theta divisors in many other situations. These applications will
be discussed elsewhere, in the thesis of the first author.

The outline of the paper is as follows. Section 1 is a general discussion of how to calculate
the multiplicity of points of $\operatorname{Sing}\Xi$. The basic idea is to use test curves to
determine the multiplicity, and to measure the multiplicity of $\Xi$ along such a curve
by relating it to the extendability of sections to infinitesimal neighborhoods of a
point. In Section 2, we deal with the case of a smooth curve $C$ and use the methods of
Section 1 to show that, if $C$ is a smooth curve of genus $6$ such that the theta
divisor $\Xi$ of
$P=P(\widetilde{C}, \pi)$ has a unique singular point of multiplicity $3$, then $C$ is
a plane quintic. Roughly, the idea is as follows: if $L$ is the line bundle on
$\widetilde{C}$ corresponding to the singular point, then
$L$ is fixed by the involution $\tau$ on $\widetilde{C}$ corresponding to the double cover, and
hence $L=\pi^*M$ for some line bundle $M$ on $C$, necessarily a theta characteristic. Moreover,
Clifford's theorem implies that $h^0(L) \leq 6$, with equality only if $C$ is hyperelliptic (in
which case the singular locus of $P$ is at least one-dimensional). The proof then involves
analyzing the cases $h^0(L) =2$ and $h^0(L)=4$ and showing that a unique singular point 
$\Xi$ which is of multiplicity $3$  can arise only for a plane quintic.

In Section 3, we extend the analysis to  singular curves of genus $6$, using Beauville's
extension of the Prym morphism to stable curves \cite{Beau1}. The main difficulty arises from the
fact that, in this case, a line bundle on $\widetilde{C}$ fixed by the involution is not
necessarily the pullback of a line bundle from $C$. Nonetheless, such line bundles can be
classified. Somewhat tedious arguments then lead to the same conclusion as for smooth curves,
that $C$ is a plane quintic.

Finally, in Section 4, we discuss the proof of a  result of Beauville, that
the Prym variety of an admissible double cover of a plane quintic corresponding to an
odd theta characteristic is in fact the intermediate Jacobian of a cubic threefold. This
then completes the proof of Theorem~\ref{ourthm}.

It is a pleasure to thank Arnaud Beauville  for discussions about his work on
Prym varieties, plane quintics, and cubic threefolds, and Roy Smith and Robert Varley for their
encouragement and for many detailed comments, both historical and mathematical, on an earlier
version of this paper.

\section{Theta divisors}

\subsection{Preliminaries}

We begin with an elementary lemma:

\begin{lemma} Let $H$ be a hypersurface, not necessarily reduced, defined in a
neighborhood of
$0$ in
$\Cee^n$ and containing $0$. Let $f\colon \Delta \to \Cee^n$ be a holomorphic map, with
$f(0) = 0$. Then $\operatorname{mult}_0H \leq \operatorname{mult}_0f^*H$, and equality
holds if and only if $f_*(T_0\Delta)$ is not contained in the tangent cone to $H$ at
$0$.\qed
\end{lemma}

Let $C$ be a smooth projective curve of genus $g$, let $S$ be a scheme, and let
$\mathcal{L}$ a line bundle over $C\times S$, of relative degree $g-1$. More generally, it suffices to assume that
$C$ is an arbitrary projective  curve and that $\chi(C;\mathcal{L}_s) =0$ for every $s$.
Then Mumford's construction
\cite[Theorem, p.\ 46]{Mumf1} gives, locally on $S$,   a complex of locally free
$\scrO_S$-modules
$\mathcal{C}^0 \xrightarrow{d} \mathcal{C}^1$  whose
cohomology is 
$R^0\pi_2{}_*\mathcal{L}$ in dimension $0$ and $R^1\pi_2{}_*\mathcal{L}$ in dimension
$1$. Suppose that $\det d$ is not a zero divisor, and set  $\Theta_S = (\det d)$. Then
$\Theta_S$ is an effective nonzero Cartier divisor on
$S$, which is independent of the choice of the complex $\mathcal{C}^\bullet$, and 
satisfies:
\begin{enumerate}
\item[(i)] The support of $\Theta_S$ is equal to the set of $s\in S$ such that $h^0(C;
\mathcal{L}_s)\neq 0$;
\item[(ii)] If $S = \Pic^{g-1}C$ and $\mathcal{L}$ is the Poincar\'e line bundle, then
$\Theta_S$ is the usual theta divisor, with multiplicity one;
\item[(iii)] The construction is functorial: if $f\colon S'\to S$ is a morphism, and
$\mathcal{L}' = (\Id\times f)^*\mathcal{L}$, then $\Theta_{S'}=f^*\Theta_S$.
\item[(iv)] If $S$ is smooth and $\dim S =1$, then $\Theta_S = \sum_{s\in
S}\ell((R^1\pi_2{}_*\mathcal{L})_s)\cdot s$.
\end{enumerate}

Clearly then we have the following, which for simplicity we shall just state in the case
where $C$ is smooth:

\begin{lemma}\label{leq} With notation as above, suppose $S$ is a smooth curve, that
$f\colon S
\to
\Pic^{g-1}C$ is the morphism induced by the line bundle $\mathcal{L}$, and that $f(s_0)
=x\in \Pic^{g-1}C$. Then $\operatorname{mult}_x\Theta \leq \deg _{s_0}\Theta_S =
\ell((R^1\pi_2{}_*\mathcal{L})_{s_0})$.  Finally, there exists an $S$ and a line
bundle  $\mathcal{L}$ as above such that equality holds.
\qed
\end{lemma}

Let $L = \mathcal{L}|C\times \{s_0\}$. Since the map $(R^1\pi_2{}_*\mathcal{L})_{s_0}
\to H^1(C; L)$ is surjective and there exists an $S$ of dimension one such that
$\operatorname{mult}_x\Theta = \deg _{s_0}\Theta_S$, we see:

\begin{lemma}\label{geq} If $L$ is the line bundle on $C$ corresponding to $x\in
\Pic^{g-1}C$, then  $\operatorname{mult}_x\Theta \geq h^1(L) = h^0(L)$.
\qed
\end{lemma}

Our goal will eventually be to apply the above in the Prym case: Let $C$ be a smooth
curve of genus $g$, $\pi \colon \widetilde{C} \to C$ a connected \'etale double cover.
More generally, we will also consider admissible double covers in the sense of
Beauville \cite{Beau1}, and leave the necessary modifications to the reader.  Let
$P\subseteq J^{2g-2}\widetilde{C}$ be the Prym variety: a line bundle
$L$ on
$\widetilde{C}$ belongs to
$P$ if and only if $\operatorname{Norm}(L) = K_C$ and $h^0(L) \equiv 0 \mod 2$. In this
case, if $\widetilde{\Theta}$ is the canonical theta divisor of
$J^{2g-2}\widetilde{C}$, then $\widetilde{\Theta}\cap P =2\Xi$, where $\Xi$ is the
class of a principal polarization. Suppose that $S$ is a smooth curve with $s_0\in S$,
and that $\mathcal{L}$ is a line bundle over $\widetilde{C}\times S$ of relative degree
$2g-2$ and such that the induced morphism $f\colon S\to J^{2g-2}\widetilde{C}$ lies in
$P$. For $x=f(s_0)$, we then have 

\begin{lemma}\label{Prymineq} With notation as above,  $$\frac12h^0(L) \leq
\operatorname{mult}_x\Xi
\leq
\frac12\deg _{s_0}\Theta_S =
\frac12\ell((R^1\pi_2{}_*\mathcal{L})_{s_0}). $$
Moreover, there exists a choice of  $S$ and a line
bundle  $\mathcal{L}$ as above such that $\operatorname{mult}_x\Xi
=\frac12\ell((R^1\pi_2{}_*\mathcal{L})_{s_0})$. \qed
\end{lemma}

\subsection{Obstructions}
We return to the case of a general curve $C$ and theta divisor $\Theta$. Let $S$ be a
smooth curve and let $t$ be a local coordinate for $S$ centered at $s_0$ and only
vanishing there.  Fix the following notation: $C_k = C\times \Spec \Cee[t]/(t^{k+1})$,
so that
$C=C_0$, and 
$\mathcal{L}_k =
\mathcal{L}/t^{k+1}\mathcal{L}$, so that $L = \mathcal{L}_0$. Thus $\mathcal{L}_k$ is
the restriction of $\mathcal{L}$ to $C_k$, and $L$ is the
restriction of $\mathcal{L}$ to $C\times \{s_0\}$. To
calculate
$\ell((R^1\pi_2{}_*\mathcal{L})_{s_0})$, we use the following:

\begin{lemma}\label{sections} For all
$N$ sufficiently large,
$\ell(H^0(C_N; \mathcal{L}_N)$ is
independent of $N$ and is equal to $\ell((R^1\pi_2{}_*\mathcal{L})_{s_0})$.
\end{lemma}
\begin{proof} Choose $N$ so large that $t^{N+1}$ annihilates 
$(R^1\pi_2{}_*\mathcal{L})_{s_0}$. Apply $R^i\pi_2{}_*$ to the short exact sequence
$$0 \to \mathcal{L} \xrightarrow{t^{N+1}}  \mathcal{L} \to \mathcal{L}_{N} \to 0.$$
By assumption $R^0\pi_2{}_*\mathcal{L}=0$. Thus we obtain an exact sequence
$$0 = R^0\pi_2{}_*\mathcal{L}\to H^0(\mathcal{L}_{N}) \to R^1\pi_2{}_*\mathcal{L}
\xrightarrow{t^{N+1}} R^1\pi_2{}_*\mathcal{L}.$$
Since $ R^1\pi_2{}_*\mathcal{L} \xrightarrow{t^{N+1}}
R^1\pi_2{}_*\mathcal{L}$ is the zero map, the homomorphism $H^0(\mathcal{L}_{N}) \to
R^1\pi_2{}_*\mathcal{L}$ is an isomorphism.
\end{proof}

Let us be more explicit about the value of $N$ above.   There is an exact sequence
$$0 \to t \mathcal{L}_k \to \mathcal{L}_k \to L \to 0,$$
where $t\mathcal{L}_k\cong \mathcal{L}_{k-1}$, and the obvious surjection
$\mathcal{L}_k \to \mathcal{L}_{k-1}$ induces a commutative diagram
$$\begin{CD}
0 @>>> t\mathcal{L}_k @>>> \mathcal{L}_k @>>> L @>>> 0\\
@. @VVV @VVV @V{=}VV @.\\
0 @>>> t\mathcal{L}_{k-1} @>>> \mathcal{L}_{k-1} @>>> L @>>> 0.
\end{CD}$$
Let $\partial _k \colon H^0(L) \to H^1(t\mathcal{L}_k) \cong H^1(\mathcal{L}_{k-1})$ be
the coboundary. Then by naturality there is a commutative diagram:
$$\begin{CD}
H^0(L) @>{\partial _k}>> H^1(t\mathcal{L}_k)\cong H^1(\mathcal{L}_{k-1})\\
@VV{=}V @VVV\\
H^0(L) @>{\partial _{k-1}}>> H^1(t\mathcal{L}_{k-1})\cong H^1(\mathcal{L}_{k-2}).
\end{CD}$$
 An easy induction then shows:

\begin{lemma}\label{coboundary} Suppose in the above notation that $\partial_{N+1}$ is
injective for some $N$. Then, for all $n\geq N$, the natural inclusion
$t^{n-N}\mathcal{L}_n
\subseteq \mathcal{L}_n$ induces an equality $H^0(t^{n-N}\mathcal{L}_n) =
H^0(\mathcal{L}_n)$. In particular $\ell(H^0(\mathcal{L}_n)) = \ell(H^0(\mathcal{L}_N))$
for all $n\geq N$.\qed
\end{lemma}

Suppose that $s\in H^0(L)$ lifts to a section $\tilde s$ of $H^0(\mathcal{L}_{k-1})$,
and hence that
$\partial_{k-1}(s) = 0$. Of course, we are free to vary $\tilde s$ by a section of
$t\mathcal{L}_{k-1}\cong \mathcal{L}_{k-2}$. By the naturality statement above, the
image of
$\partial_k(s)$ in $H^1(\mathcal{L}_{k-2})$ is equal to $\partial_{k-1}(s) =0$. It
follows that $\partial_k(s)$ lies in the image of 
$H^1(t^{k-1}\mathcal{L}_{k-1})$ in 
$H^1(\mathcal{L}_{k-1})$, using the long exact sequence associated to 
$$0 \to t^{k-1}\mathcal{L}_{k-1} \to \mathcal{L}_{k-1} \to \mathcal{L}_{k-2} \to 0.$$
Thus, using the isomorphism $t^{k-1}\mathcal{L}_{k-1} \cong L$, there is a well-defined
element
$\gamma_k(s) \in H^1(L)/\operatorname{Im}(H^0(\mathcal{L}_{k-2}))$ which represents the
obstruction to lifting $s$ to a section of $\mathcal{L}_k$. Clearly, the image of
$H^0(\mathcal{L}_{k-2})$ represents the freedom to alter $\tilde s$ by a section of
$t\mathcal{L}_{k-1}\cong \mathcal{L}_{k-2}$. In particular, if $\partial _1 = \cdots =
\partial _{k-1}=0$, then $\gamma_k(s)$ is a well-defined element of $H^1(L)$.

For $k=1$, we set $\gamma_1(s) =\partial
_1(s)$. In this case, a well-known argument gives:

\begin{lemma}\label{cupprod} Let $\xi \in H^1(\scrO_C) = \Ext^1(L, L)$ be the extension
class corresponding to $\mathcal{L}_1$. Then $\partial_1(s) = s\cup \xi$, where the cup
product is $H^0(L) \otimes H^1(\scrO_C) \to H^1(L)$. \qed
\end{lemma}

For computation, we will use the following which is an immediate consequence of the
usual properties of cup product:

\begin{lemma}\label{cupprod2} Let $D$ be an effective divisor on $C$ and let $\partial
\colon H^0(\scrO_D(D)) \to H^1(\scrO_C)$ be the coboundary map induced by the short
exact sequence
$$0 \to \scrO_C \to \scrO_C(D) \to \scrO_D(D) \to 0.$$
Suppose that $\xi \in H^1(\scrO_C)$ is of the form $\partial(t)$ for some $t\in
H^0(\scrO_D(D))$. Then $s\cup \xi = \partial_L (s\cdot t)$, where $s\cdot t$ is the
section of $L(D)|D$ given by taking the cup product of $s$ and $t$, and $\partial_L$
is the coboundary  homomorphism arising from
$$0 \to L \to L\otimes \scrO_C(D) \to L\otimes \scrO_D(D) \to 0.\qed$$
\end{lemma}

As a consequence:

\begin{theorem}[Riemann Singularity Theorem] If $L$ is the line bundle on $C$
corresponding to $x\in J^{g-1}C$, then 
$\operatorname{mult}_x\Theta  = h^0(L)$.
\end{theorem}
\begin{proof} By Lemma~\ref{geq}, $\operatorname{mult}_x\Theta  \geq h^0(L)$. To see 
the opposite inequality, suppose that we can find a $\xi \in H^1(\scrO_C)$ such that
$\cdot \cup\xi \colon H^0(L) \to H^1(L)$ is injective. By Lemma~\ref{cupprod} and
Lemma~\ref{coboundary}, it follows that, given $S$ smooth of dimension one,  $s_0 \in
S$, and  $f\colon S \to \Pic^{g-1}C$ such that $f(s_0) = x$ and with tangent vector
$\xi$, $\deg (\Theta_S)_{s_0} = h^0(L)$. Since by Lemma~\ref{leq} 
$\operatorname{mult}_x\Theta \leq \deg (\Theta_S)_{s_0} = h^0(L)$, we conclude that
$\operatorname{mult}_x\Theta = h^0(L)$.

To find the appropriate $\xi$, apply Lemma~\ref{cupprod2} with $D=\sum_{i=1}^{r+1}p_i$,
where $h^1(L) = r+1$ and the $p_i$ are general points of $C$, and
$t$ is a section of $\scrO_D(D)$ which does not vanish at any of the $p_i$. If the
$p_i$ are general, then $H^0(L\otimes \scrO_C(-D)) = 0$. In particular $s\cdot t$ is
nonzero for every nonzero
$s\in H^0(L)$. It suffices to show that $\partial _L \colon H^0(L\otimes \scrO_C(D)|D)
\to H^1(L)$ is injective, or equivalently that $H^0(L) \to H^0(L\otimes \scrO_C(D))$ is
surjective, i.e.\ that
$\dim H^0(L\otimes \scrO_C(D)) =r+1$. By Riemann-Roch it suffices to show that
$H^1(L\otimes \scrO_C(D)) =0$ or equivalently that $H^0(K_C\otimes L^{-1}\otimes
\scrO_C(-D))=0$. Since $h^0(K_C\otimes L^{-1}) = r+1$, if the $p_i$ are general then
indeed $H^0(K_C\otimes L^{-1}\otimes
\scrO_C(-D))=0$ and we are done.
\end{proof}

\subsection{Cocycle calculations}
We conclude with some computations. For an appropriate affine cover $\{U_i\}$ of $C$,
we may assume that $L$ has transition functions $\lambda_{ij}$ and that the transition
functions for $\mathcal{L}$ are of the form $\lambda_{ij}(t) = \lambda_{ij}(1 +
\sum_{k=1}^\infty\alpha_{ij}^{(k)}t^k)$. Here $\xi = \alpha_{ij}^{(1)}$ is a
$1$-cocycle for $\scrO_C$. Likewise set 
$$\lambda_{ij;N}(t) = \lambda_{ij}(1 +
\sum_{k=1}^N\alpha_{ij}^{(k)}t^k),$$
a $1$-cocycle for $\scrO_{C_{N}}^*$. 
Assume that $s\in H^0(L)$ and that $s_{N-1}$ is
a lifting of
$s$ to a section of $\mathcal{L}_{N-1}$. Thus using the trivialization over the
open cover $\{U_i\}$ we have  
$$s_{i;N-1} = s_i + \sum _{i=1}^{N-1}\sigma_i^{(k)}t^k$$
for some functions $\sigma_i^{(k)}$ on $U_i$,
with $s_{i;N-1} =\lambda_{ij;N-1}(t)s_{j;N-1}$ on $(U_i\cap U_j) \times
\Spec\Cee[t]/(t^{N})$. Set $\sigma_i^{(0)} =s_i$. The section $s_{N-1}$ lifts to a
section
$s_N$ if and only if there exists $\sigma_i^{(N)}\in \Gamma \scrO_{U_i}$ such that, if
we set $s_{i;N} = s_i +
\sum _{i=1}^{N}\sigma_i^{(k)}t^k$, then $s_{i;N} =\lambda_{ij;N}(t)s_{j;N}$ on
$(U_i\cap U_j) \times
\Spec\Cee[t]/(t^{N+1})$. Since $s_{i;N-1}$ is already a section of
$\mathcal{L}_{N-1}$, this becomes the equality
$$\sigma_i^{(N)} = \lambda_{ij}\sigma_j^{(N)} + \sum
_{k=0}^{N-1}\lambda_{ij}\alpha_{ij}^{(N-k)}\sigma_j^{(k)}.$$
In other words,

\begin{lemma}\label{calcobstr} The section $s_{N-1}$ lifts to a section $s_N$ if and
only if the image of $\sum
_{k=0}^{N-1}\lambda_{ij}\alpha_{ij}^{(N-k)}\sigma_j^{(k)}$ in $H^1(L)$ is zero. More
generally, this image is the obstruction $\gamma_k(s)$ defined above.
\qed
\end{lemma}

Suppose that $p$ is a point of $C$, and fix once and for all a local coordinate $z$ at
$p$. More precisely, let $\{U_i\}$ be an open cover of $C$, and assume that $p\in U_0$,
that $p\notin U_i$ for $i\neq 0$, and that $z\in \Gamma\scrO_{U_0}$ is a coordinate
centered at $p$. A calculation then shows:

\begin{lemma}\label{1.10} For $a\in \Cee$, let $\xi\in H^1(\scrO_C)$ be the image of
$a/z$ under the coboundary map induced by the short exact sequence
$$0\to \scrO_C \to \scrO_C(p) \to \scrO_C(p)|p \to 0.$$
Let $s\in H^0(L)$ be a section such that $s(p) =0$. Then $s\cup \xi =0$ in $H^1(L)$,
and in fact, choosing $\xi_{ij}$ to be given by the $1$-cocycle
$$
\xi_{ij} = \begin{cases}a/z, &\text{if $i=0$;}\\ 
0, &\text{otherwise,}
\end{cases}$$
then $s\cup \xi = \delta\sigma$, where $\sigma$ is the $0$-cochain defined by
$\sigma_0 = as/z$ and $\sigma_j=0$ for $j\neq 0$, and $\delta$ is the \v Cech coboundary
map.\qed
\end{lemma}

We turn now to obstructions of order $2$. More generally, consider the case where
$S=C$, or an affine open subset of $C$ containing the point $p$, and $\mathcal{L}$ is
either the line bundle $\pi_1^*(L\otimes \scrO_C(p))\otimes  \scrO_{C\times C}(-\Delta)$
or the line bundle $\pi_1^*(L\otimes \scrO_C(-p))\otimes \scrO_{C\times C}(\Delta)$,
where $\Delta \subseteq C\times C$ is the diagonal. We fix the coordinate $z$ as before,
and let
$t$ be the coordinate
$z$, viewed as a coordinate on an affine open subset of $C$. In the first case,
transition functions for
$\mathcal{L}$ are of the form $\lambda_{ij}$, if neither $i$ nor $j$ is zero, whereas
(for small $t$),
$$\lambda_{0j}(t) = \lambda_{0j} \cdot \left(\frac{z}{z-t}\right)= \lambda_{0j}
\cdot\sum_{k=0}^\infty \left(\frac{t^k}{z^k}\right) .$$
 In the second case,
$\lambda_{ij}(t)=\lambda_{ij}$  if neither
$i$ nor
$j$ is zero, and (for small $t$),
$$\lambda_{0j}(t) = \lambda_{0j} \cdot
\left(\frac{z-t}{z}\right)= \lambda_{0j} \cdot
\left(1 - \left(\frac{t}{z}\right)\right).$$
In case $\mathcal{L} = \pi_1^*(L\otimes \scrO_C(p))\otimes \scrO_{C\times C}(-\Delta)$,
a first order lift corresponds to $\sigma_i^{(1)}$ satisfying
$$\sigma_0^{(1)} = \lambda_{0j}\sigma_j^{(1)} +
\lambda_{0j}\left(\frac{1}{z}\right)s_j = \lambda_{0j}\sigma_j^{(1)} + \frac{s_0}{z},$$
and the second order obstruction satisfies
$$\gamma_2(s)_{0j} = \lambda_{0j}\left(\frac{1}{z^2}\right) s_j
+\lambda_{0j}\left(\frac{1}{z}\right)\sigma_j^{(1)} =
\frac{s_0}{z^2} +\left(\frac{\sigma_0^{(1)}}{z} -
\frac{s_0}{z^2}\right) = \frac{\sigma_0^{(1)}}{z}.$$
Likewise, in case $\mathcal{L} = \pi_1^*(L\otimes \scrO_C(-p))\otimes \scrO_{C\times
C}(\Delta)$, a first order lift corresponds to $\sigma_i^{(1)}$ satisfying
$$\sigma_0^{(1)} = \lambda_{0j}\sigma_j^{(1)} +
\lambda_{0j}\left(-\frac{1}{z}\right)s_j = \lambda_{0j}\sigma_j^{(1)} - \frac{s_0}{z},$$
and the second order obstruction satisfies
$$\gamma_2(s)_{0j} = \lambda_{0j}\cdot 0\cdot s_j
+\lambda_{0j}\left(-\frac{1}{z}\right)\sigma_j^{(1)} = -\frac{\sigma_0^{(1)}}{z} -
\frac{s_0}{z^2}.$$

\begin{lemma}\label{1.11} With notation as above, suppose that $s\in H^0(L)$ lifts to a
section of
$\mathcal{L}_1$, corresponding to the $0$-cochain $\sigma_i^{(1)}$.
Then $\gamma_2(s)= \{\gamma_2(s)_{ij}\}$ satisfies: $\gamma_2(s)_{ij} =0$ if neither
$i$ nor $j$ is zero, and
$$\gamma_2(s)_{0j} =
\begin{cases}
\sigma_0^{(1)}/z, &\text{if
$\mathcal{L} = \pi_1^*(L\otimes \scrO_C(p))\otimes \scrO_{C\times C}(-\Delta)$;}\\
 -\sigma_0^{(1)}/z - s_0/z^2 &\text{if $\mathcal{L}
=\pi_1^*(L\otimes
\scrO_C(-p))\otimes \scrO_{C\times C}(\Delta)$.}
\end{cases}$$
\end{lemma}
\begin{proof} This follows  from Lemma~\ref{calcobstr} and the above
calculations.
\end{proof}

An immediate calculation then gives:

\begin{corollary}\label{1.12} Let $\partial =\partial_{L, 2p}$ be the coboundary map
$$\partial \colon H^0(L\otimes \scrO_C(2p)|2p)
\to H^1(L)$$ induced from the exact sequence
$$0 \to L \to L\otimes \scrO_C(2p) \to L\otimes \scrO_C(2p)|2p\to 0,$$
and let $a\in H^0(L\otimes \scrO_C(2p)|2p)$ be defined by
$$a = 
\begin{cases}
\sigma_0^{(1)}/z,   &\text{if
$\mathcal{L} =
\pi_1^*(L\otimes
\scrO_C(p))\otimes \scrO_{C\times C}(-\Delta)$;}\\ 
-\sigma_0^{(1)}/z -
s_0/z^2, &\text{if $\mathcal{L}
=\pi_1^*(L\otimes
\scrO_C(-p))\otimes \scrO_{C\times C}(\Delta)$.}
\end{cases}$$
Then $\gamma_2(s)$ is equal to the image of $\partial (a)$ in $H^1(L)/\partial _1(  H^0(L))$.
\qed
\end{corollary}

\begin{remark} We shall use the following two observations:
\begin{enumerate}
\item[(i)] The same calculations work for singular (reduced) curves $C$, provided that
we use a smooth point $p\in C$.
\item[(ii)] The above calculations are local, in the sense that if 
$$\mathcal{L} =
\scrO_C(\pm(p_1-q_1)) \otimes \cdots \otimes \scrO_C(\pm(p_k-q_k)),$$
 where $q_i=q_i(s)$
for some $s\in S$ with $q_i(s_0) = p_i$, then the  obstructions are sums of local
contributions as calculated in the preceding lemmas.
\end{enumerate}
\end{remark}

\section{The Prym case}

We keep the following notation:   $C$ is a smooth
curve of genus $g$ with connected \'etale double
cover $\pi \colon \widetilde{C} \to C$ and involution $\tau\colon
\widetilde{C} \to \widetilde{C}$, and   $\eta$ is the corresponding line bundle of order
$2$ on $C$. The Prym variety is denoted $P$ and its canonical theta divisor is $\Xi$.

\subsection{The balanced case}

We begin with a general criterion, along the lines of the Riemann singularity theorem,
for when $\operatorname{mult}_x\Xi = n$.

\begin{lemma}\label{21}  Let $L$ be a line bundle on $\widetilde{C}$ corresponding
to the point 
$x\in P$, and suppose that $h^0(L) = 2n$. Then $\operatorname{mult}_x\Xi \geq
\frac12h^0(\widetilde{C}; L)=n$. Moreover, if there exist distinct points $q_1, \dots,
q_n \in C$, with $\pi^{-1}(q_i) =\{p_i, \tau(p_i)\}$, such that the $2n$ points   $p_1,
\tau(p_1), \dots, p_n, \tau(p_n)$ impose independent conditions on $L$, then
$\operatorname{mult}_x\Xi = n$.
\end{lemma}
\begin{proof} By Lemma~\ref{Prymineq}, $\operatorname{mult}_x\Xi \geq 
\frac12h^0(\widetilde{C}; L)$. Moreover,  equality holds if and only if  
$\operatorname{mult}_x\Xi \leq n$, if and only if there exists $\xi\in
H^1(\scrO_{\widetilde{C}})^-$ such that the cup product map $\cup\xi\colon
H^0(\widetilde{C}; L)\to H^1(\widetilde{C}; L)$ is injective.

Let $q_1, \dots, q_k$ be  points of $C$ and let
$\pi^{-1}(q_i) = \{p_i, \tau(p_i)\}$. Set $D =\sum _{i=1}^k(p_i+\tau(p_i))$ and let 
$t$ be a section of $\scrO_{\widetilde{C}}(D)|D$ which is nonzero at each point in the
support of $D$ and such that $\tau(t) =-t$.  Let $\xi=\partial (t)$, so that $\tau^*\xi
=-\xi$. The map
$\cup \xi \colon H^0(L) \to H^1(L)$ is the composition of the maps $H^0(L) \to
H^0(L\otimes \scrO_{\widetilde{C}}(D)|D)$ given by multiplying by $t$, and the
coboundary $\partial_{L,D}\colon H^0(L\otimes \scrO_{\widetilde{C}}(D)|D) \to
H^1(L)$ arising from the long exact sequence
$$0 \to L \to L\otimes \scrO_{\widetilde{C}}(D) \to (L\otimes
\scrO_{\widetilde{C}}(D))|D \to 0.$$
Note that the map $s\in H^0(L) \mapsto st$ is injective if and only if
$H^0(L(-D))=0$, and $\partial_{L,D}$ is injective if and only if $H^0(L) =
H^0(L(D))=2n$. By Riemann-Roch, $h^0(L(D)) = h^1(L(D)) + 2k$. Thus, if $\cup \xi$ fails
to be injective, then either $h^0(L\otimes
\scrO_{\widetilde{C}}(-D))\neq 0$ or $h^1(L\otimes
\scrO_{\widetilde{C}}(D))> 2n-2k$. Since  by Serre duality $h^1(L\otimes
\scrO_{\widetilde{C}}(D)) = h^0(\tau^*L\otimes
\scrO_{\widetilde{C}}(-D)) = h^0(L\otimes
\scrO_{\widetilde{C}}(-D))$, for $k\leq n$ this last condition is equivalent to the
condition that the $2k$ points $p_1, \dots, p_k, \tau(p_1),
\dots, \tau(p_k)$ fail to impose independent conditions on the free part of $L$.
Likewise for $k=n$,  
$H^0(L\otimes
\scrO_{\widetilde{C}}(-D)) \neq 0$ if and only if the $2n$ points $p_1, \dots, p_n,
\tau(p_1),
\dots, \tau(p_n)$ fail to impose independent conditions on the free part of $L$.
Hence, if the points   $p_1,
\tau(p_1), \dots, p_n, \tau(p_n)$ impose independent conditions on $L$, then $\cup \xi$
is injective and hence $\operatorname{mult}_x\Xi  =n$.
\end{proof}

\begin{remark} An elementary argument shows the following: if
$h^0(L) = 2n$ but, for every choice of distinct points $q_1,
\dots, q_n \in C$ with $\pi^{-1}(q_i) =\{p_i, \tau(p_i)\}$,  the $2n$
points   $p_1,
\tau(p_1), \dots, p_n, \tau(p_n)$ fail to impose independent conditions on $L$, then
$L$ is of the form $\pi^*M(E)$, where $M$ is a line bundle on $C$ with $h^0(M) \geq 2$
and $E$ is effective.
\end{remark}

\begin{theorem}[Smith-Varley \cite{SV2}]\label{SV} Let $L$ be a line bundle on
$\widetilde{C}$ corresponding to the point 
$x\in P$. Suppose that $L=\pi^*M\otimes
\scrO_{\widetilde{C}}(B)$, where $h^0(L) = h^0(\pi^*M)$, so that $B$ is contained in
the base locus of $L$. If  $h^0(M) = h^0(M\otimes
\eta)$, then  $\operatorname{mult}_x\Xi = n$. If $h^0(M) \neq h^0(M\otimes
\eta)$, then, for every line bundle $\mathcal{L}_1$  over $C\times
\Spec\Cee[t]/(t^2)$ restricting to $L$ on the closed fiber and corresponding to a
morphism $\Spec\Cee[t]/(t^2) \to P$, $\ell(h^0(\mathcal{L}_1)) \geq 2\max (h^0(M) ,
h^0(M\otimes \eta))$. Hence
$$\operatorname{mult}_x\Xi \geq \max (h^0(M) , h^0(M\otimes \eta))>\frac12h^0(L).$$
\end{theorem}
\begin{proof} 
By assumption, $L=\pi^*M\otimes
\scrO_{\widetilde{C}}(B)$, where $h^0(L) = h^0(\pi^*M)$. First suppose that $h^0(M) =
h^0(M\otimes\eta) =n$. Note that 
$$K_{\widetilde{C}}\otimes
L^{-1}= \tau^*L =\pi^*M\otimes \scrO_{\widetilde{C}}(\tau(B)).$$
Thus $h^0(K_{\widetilde{C}}\otimes
L^{-1}) = h^0(\pi^*M)$ as well.    Choose $n$ general points $q_1,
\dots, q_n$  on $C$ and write
$\pi^{-1}(q_i) = \{p_i, \tau(p_i)\}$. Let $D =\sum _{i=1}^n(p_i+\tau(p_i))$, let $t$
be a section of $\scrO_{\widetilde{C}}(D)|D$ which is nonzero at each point in the
support of $D$ and such that $\tau(t) =-t$, and let $\partial(t) =\xi$.  It suffices to
find $q_i$ so that 
$$H^0(L\otimes
\scrO_{\widetilde{C}}(-D)) =H^1(L\otimes
\scrO_{\widetilde{C}}(D)) = 0.$$
We have seen in the course of the proof of Lemma~\ref{21} that both groups have the
same dimension. Working for example with $H^0(L\otimes
\scrO_{\widetilde{C}}(-D))$, we have
\begin{align*}h^0(L\otimes \scrO_{\widetilde{C}}(-D)) &= h^0(\pi^*(M\otimes
\scrO_C(-\sum_iq_i)))\\ &=
h^0(M\otimes\scrO_C(-\sum_iq_i)) + h^0(M\otimes\eta\otimes\scrO_C(-\sum_iq_i)) .
\end{align*}
Since $h^0(M) = h^0(M\otimes\eta) =n$, if the $q_i$ are general then both of the
terms in the last line are zero as claimed. Hence $H^0(L\otimes
\scrO_{\widetilde{C}}(-D)) = H^1(L\otimes
\scrO_{\widetilde{C}}(D)) = 0$.

Finally assume that $L=\pi^*M(B)$, with $h^0(M) = n_1 > h^0(M\otimes \eta) = n_2$.
 For $\xi \in H^1(\scrO_{\widetilde{C}})^-$, consider
the pairing $H^0(L) \otimes H^0(\tau^*L) \to \Cee$ defined by
$$\langle s, t \rangle = (st)\cup \xi \in H^1(L\otimes \tau^*L)= H^1(K_{\widetilde{C}})
\cong
\Cee.$$ Since $\tau$ acts trivially on $H^1(K_{\widetilde{C}})$, we must have $\langle
\tau(t), \tau(s)\rangle = -\langle s, t \rangle$. By Serre duality, it suffices to show
that
$\langle \cdot , \cdot \rangle$ is degenerate in the first variable for every $\xi$.
But $H^0(L) = H^0(\pi^*M)$ and  $H^0(\tau^*L)= H^0(\pi^*M(\tau(B)))= H^0(\pi^*M)$, so
that $\langle \cdot , \cdot \rangle$ induces a  pairing $H^0(\pi^*M) \otimes
H^0(\pi^*M)
\to \Cee$.   Clearly, given sections $s,t\in H^0(\pi^*M)$, 
$\langle i(s), i(t)
\rangle = i(st)\cup \xi$, where we denote by $i$ any of the natural inclusions
$\pi^*M\to \pi^*M(B)$, $\pi^*M\to \pi^*M(\tau(B))$, or $\pi^*M^{\otimes 2} \to
\pi^*M(B) \otimes \pi^*M(\tau(B))$, and thus $\langle i(s), i(t) \rangle =  \langle
i(t),i( s) \rangle$ for all $s,t\in H^0(\pi^*M)$. Hence, for $s,t\in H^0(\pi^*M) \otimes
H^0(\pi^*M)$, $\langle
\tau(i(s)), \tau(i(t))\rangle = -\langle i(s), i(t) \rangle$. It follows that $\langle
\cdot ,
\cdot \rangle$ pairs $H^0(M)$ with $H^0(M\otimes \eta)$ and in fact has a kernel of
dimension at least $n_1-n_2$. Hence $\cup \xi$ has a kernel of dimension at least
$n_1-n_2$. Thus, if $\mathcal{L}_1$ is the line bundle over $C\times
\Spec\Cee[t]/(t^2)$ corresponding to $\xi$, $\ell(\mathcal{L}_1) \geq
(n_1+n_2)+(n_1-n_2) =2n_1$. The final statement of Theorem~\ref{SV} is then clear.
\end{proof}

\begin{remark} As pointed out to us by Smith and Varley, the argument above shows that, if $L =
\pi^*M(B)$ for an effective divisor $B$ on $\widetilde{C}$, then the multiplicity of $\Xi$ at $L$
is at least
$h^0(M)$. Hence, if
$h^0(M) > \frac12h^0(L)=n$, then the multiplicity of $\Xi$ at $L$ is greater than $n$. In fact,
Smith and Varley show in \cite{SV2} that, if $h^0(L) = 2n$, then the multiplicity of
$\Xi$ at $L$ is greater than $n$ if and only if there exists a line bundle $M$ on $C$ with
$h^0(M) > n$ and an effective divisor $B$ on $\widetilde{C}$ such that $L =
\pi^*M(B)$.
\end{remark}

\subsection{The case where all sections are invariant}
We consider now the opposite extreme, where all of $H^0(L)$ is invariant. Our goal
is to show the following, special cases of which were established by Smith-Varley:

\begin{theorem}\label{h0istwo} Suppose that $L=\pi^*M(B)$ corresponds to the point $x\in
P$ and that $h^0(L) = h^0(M)=2n>0$. Then
$\operatorname{mult}_x\Xi = 2n$.
\end{theorem}
\begin{proof} Note that the hypothesis implies (in fact is equivalent to) the statement that $B$
is contained in the base locus of
$L$, and that, if
$p\in \widetilde{C}-B$ and
$s$ is a section of $L$ vanishing at $p$, then $s$ vanishes at $\tau(p)$ also. We begin by
showing that every section of
$L$ lifts to first order:

\begin{lemma} Suppose that $L=\pi^*M(B)$ is as above  and
that $h^0(L) = h^0(M)$. Let $\mathcal{L}_1$ be a line bundle over $C\times
\Spec\Cee[t]/(t^2)$ restricting to $L$ on the closed fiber and corresponding to a
morphism $\Spec\Cee[t]/(t^2) \to P$. Then the homomorphism $H^0(\mathcal{L}_1) \to
H^0(L)$ is surjective. Hence $\ell(H^0(\mathcal{L}_1)) = 2h^0(L)$ and
$\operatorname{mult}_x\Xi \geq 2n$.
\end{lemma}
\begin{proof} By Theorem~\ref{SV} and its proof, the homomorphism $\cup\xi\colon H^0(L)
\to H^1(L)$ is zero for every $\xi\in H^1(\scrO_{\widetilde{C}})^-$. Thus, if $\mathcal{L}_1$
is the line bundle over $C\times
\Spec\Cee[t]/(t^2)$ corresponding to $\xi$, then  
$H^0(\mathcal{L}_1) \to
H^0(L)$ is surjective  and $\ell(H^0(\mathcal{L}_1)) = 2h^0(L)$. The last statement is
clear.
\end{proof}

In particular, for every choice of $\mathcal{L}_2$, $\partial_1=0$ and so $\gamma_2(s)$ is a
well-defined element of $H^1(L)$, independent of the choice of a lift of $s$ to a section of
$\mathcal{L}_1$. 

\begin{lemma}\label{newlemma}  Let $h^0(L) = h^0(M) = 2n$. Suppose that $q_1, \dots, q_n$ are
general points of
$C$, with $\pi^{-1}(q_i) = \{p_i, \tau(p_i)\}$. Let $D = \sum _i(p_i+\tau(p_i))$. Then the
natural inclusions 
$$H^0(L) \subseteq H^0(L(\sum_{i=1}^np_i) )\qquad \textrm{and} \qquad H^0(L(D)) \subseteq
H^0(L(2D-\sum_{i=1}^np_i))$$ are equalities. A similar statement holds with $p_i$ replaced by
$\tau(p_i)$.
\end{lemma}
\begin{proof} Clearly, the $q_i$ are general points of $C$ if and only if the $p_i$ are general
points of $\widetilde{C}$. To see the equality
$H^0(L) = H^0(L(\sum_{i=1}^np_i) )$, note that Riemann-Roch implies that
$$h^0(L(\sum_{i=1}^np_i) ) - h^1(L(\sum_{i=1}^np_i) ) = n.$$
But $h^1(L(\sum_{i=1}^np_i) ) = h^0(\tau^*L(-\sum_{i=1}^np_i) )$. Since the $n$ points $p_i$ are
general, they impose independent conditions on the linear series $|\tau^*L|$, and since
$h^0(\tau^*L)=h^0(L) = 2n$, $h^0(\tau^*L(-\sum_{i=1}^np_i) ) =n$. It follows that
$h^0(L(\sum_{i=1}^np_i) ) =2n = h^0(L)$, and so the inclusion $H^0(L) \subseteq
H^0(L(\sum_{i=1}^np_i) )$ is an equality.

The equality $H^0(L(D)) \subseteq
H^0(L(2D-\sum_{i=1}^np_i))$ is similar: by Riemann-Roch,
$$h^0(L(D)) = 2n + h^1(L(D)) = 2n + h^0(\tau^*L(-D)).$$
If $D$ does not meet $B$, then 
$$h^0(\tau^*L(-D)) = h^0(\pi^*M(-D)) = h^0(M(-\sum_iq_i)) = n$$ 
if the $q_i$ are general, since $h^0(M) = 2n$. Thus $h^0(L(D)) = 3n$. Likewise,
\begin{align*}
h^0(L(2D-\sum_ip_i))& = 3n + h^1(L(2D-\sum_ip_i))\\ = 3n +
h^0(\tau^*L(-2D+\sum_ip_i))
&= 3n+h^0(\pi^*M(-\sum_i(p_i+2\tau(p_i)))).
\end{align*}
Since $H^0(\pi^*M)=H^0(M)$, $h^0(\pi^*M(-\sum_i(p_i+2\tau(p_i)))) = h^0(M(-2\sum_iq_i))$.
Clearly, for $n$ general points $q_i$, $h^0(M(-2\sum_iq_i)) =0$. Thus 
$h^0(L(2D-\sum_ip_i)) = 3n =h^0(L(D))$ and we conclude as before.
\end{proof}

We now complete the proof of Theorem~\ref{h0istwo}. Choose $n$ general points  $q_i$ of
$C$ and set $\pi^{-1}(q_i) = \{p_i, \tau (p_i)\}$ and $D = \sum_{i=1}^n (p_i +
\tau(p_i))$.  Let $\Delta$ denote a disk in $\Cee$ and let
$\mathcal{L}$ be the line bundle over $\widetilde{C}\times \Delta$ whose restriction to the
fiber over $t$ is  $L\otimes
\scrO_{\widetilde{C}}\left(\sum _{i=1}^n(p_i-\tau(p_i)-r_i(t) + \tau(r_i(t)))\right)$, where the
$r_i\colon
\Delta \to \widetilde{C}$ are isomorphisms from $\Delta$ to a neighborhood of $p_i$ such that
$r_i(0)=p_i$.  Let $z_i$ be the local $\tau$-invariant coordinate  around $p_i$ and
$\tau(p_i)$ given by $r_i^{-1}$. Choose a sufficiently small open cover $\{U_\alpha\}_{\alpha\in
I}$ of
$\widetilde{C}$ such that each $U_\alpha$ contains at most one of the points $p_i, \tau(p_i)$.
Given $i$, let $\alpha_i$ be the unique element of $I$ such that 
$p_i\in U_{\alpha_i}$ and  let  $\beta_i$ be the unique element of $I$ such that 
$\tau(p_i)\in U_{\beta_i}$. We may assume that $U_{\alpha_i}\cap U_{\alpha_j} = U_{\beta_i}\cap
U_{\beta_j} =\emptyset$ for $i\neq j$, and also that $U_{\alpha_i}\cap U_{\beta_j} =
\emptyset$ for all $i,j$.  Then
$\mathcal{L}_1$ is the line bundle defined by
$\xi\in H^1(\scrO_{\widetilde{C}})^-$, where $\xi_{\alpha_i \gamma} = 1/z_i$, 
$\xi_{\beta_i\gamma } = -1/z_i$, and $\xi_{\mu\nu}$ is
otherwise
$0$.  We have seen that every section
$s$ of
$L$ lifts to first order. Thus it suffices to show that $\gamma_2(s)\neq 0$ for every
nonzero section
$s$. Here $\gamma_2(s)$ is a well-defined element of $H^1(L)$, given by
Lemma~\ref{1.11}. By Corollary~\ref{1.12}, if $\gamma_2(s) = 0$, then $\partial_{L,2D}
(a) = 0$, where $a\in H^0(L\otimes \scrO_{\widetilde{C}}(2D)|2D)$ given
as follows: if  $s\cup \xi =\delta
(\sigma^{(1)})$, then
$$a_{p_i} = \sigma^{(1)}_{\alpha_i}/z_i; \quad a_{\tau(p_i)}
= -\sigma^{(1)}_{\beta_i}/z_i  -
s/z_i^2.$$
If $\partial_{L,2D}(a) = 0$, then there exists a section $\Sigma\in H^0(L(2D))$ whose image in
$H^0(L(2D)|2D)$ is $a$. It follows from the explicit description of $a$ that, in fact, $\Sigma
\in H^0(L(2D-\sum_ip_i))$. By Lemma~\ref{newlemma}, $\Sigma \in H^0(L(D))$, in other words
$s(\tau(p_i)) =0$ for every $i$. By the remarks at the begininng of the proof, $s(p_i)=0$
for every $i$ as well. By Lemma~\ref{1.10}, we can then choose
$\sigma^{(1)}$ as follows:  $\sigma^{(1)}_{\alpha_i} = s/z_i$,  
$\sigma^{(1)}_{\beta_i} = -s/z_i$, and otherwise $\sigma^{(1)}_\gamma=0$. Thus $a_{\tau(p_i)} =
s/z_i^2-s/z_i^2=0$, and so $\Sigma$ has no pole at $\tau(p_i)$. But then $\Sigma \in
H^0(L(\sum_ip_i))= H^0(L)$, by Lemma~\ref{newlemma} again. It then follows that
$\sigma^{(1)}_{\alpha_i}/z_i = s/z_i^2=0$ as a section of $L(2D)|2p_i$ for every $i$, in
other words that
$s$ vanishes to order $2$ at $p_i$. Using $H^0(L) =H^0(M)$ and the fact that the
$p_i,\tau(p_i)$ are not in $B$, $s$ also vanishes to order $2$ at $\tau(p_i)$, so that
$s\in H^0(L(-2D))$. Finally, we have already seen that $H^0(L(-2D)) = H^0(M(-2\sum_iq_i)) =0$.
Thus, if $\gamma_2(s) =0$, then $s=0$, so that $\gamma_2$ is injective.  This completes the proof
of Theorem~\ref{h0istwo}.
\end{proof}

Given Theorem~\ref{h0istwo} and the theorem of Smith-Varley, it is natural to make the following
conjecture (which has also been made by Smith-Varley):

\begin{conjecture} Suppose that $\pi \colon \widetilde{C}\to C$ is a connected \'etale double
cover of the smooth curve $C$ and that $L$ is a line bundle corresponding to the point $x\in
P(\widetilde{C}, \pi)$ with $h^0(L) = 2n$ and $\operatorname{mult}_x\Xi > n$. Let $M$ be a line
bundle on $C$ such that $L=\pi^*M(B)$, with $h^0(M)$ maximal. Then $\operatorname{mult}_x\Xi
=h^0(M)$.
\end{conjecture}

We can now prove the main theorem for Prym varieties of smooth curves of genus $6$. 

\begin{theorem}\label{smoothcase} Let
$\pi \colon \widetilde{C}\to C$ be a connected \'etale double cover of a smooth curve of
genus
$6$. If the theta divisor of the associated Prym variety $P(\widetilde{C}, \pi)$ has a
unique singular point which is an isolated triple point, then $C$ is a plane quintic
and $\pi\colon \widetilde{C}\to C$ is the double cover associated to an odd theta
characteristic.
\end{theorem}
\begin{proof} Let $L$ be a line bundle on $\widetilde{C}$ of degree $10$ corresponding
to a singular point of the theta divisor $\Xi$ of $P(\widetilde{C}, \pi)$. In
particular,
$h^0(L) = 2n > 0$. By Clifford's theorem, $h^0(L)\leq 6$, with equality if and only if
$\widetilde{C}$ is hyperelliptic. This easily implies that $C$ is hyperelliptic (see
below for the singular case), and then $\Xi$ does not have an
isolated singularity. Thus we may assume that $h^0(L)< 6$, and hence $h^0(L) \leq 4$.

Since there is a unique singular point of $\Xi$, $\tau^*L = L$, and by \'etale descent
$L=\pi^*M$. Thus $h^0(L) = h^0(M) + h^0(M\otimes \eta)\leq 4$, where $\eta$ is the line
bundle of order two corresponding to $\widetilde{C}$. We may assume that 
$h^0(M)
\geq h^0(M\otimes \eta)$ and that $h^0(M)\neq 0$. By Clifford's theorem, $h^0(M) \leq
5/2 + 1 < 4$. If 
$h^0(M) = 1$, then the only possibility is  $h^0(M\otimes \eta) =1$ as well. In this
case the multiplicity of the corresponding point of $\Xi$ is
$h^0(M) =1$, by  Theorem~\ref{SV}. If  $h^0(M) = 2= h^0(M\otimes \eta)$, then the 
multiplicity of the corresponding point of $\Xi$ is
$h^0(M) =2$, again by
Theorem~\ref{SV}. If $h^0(M) = 2$ and
$h^0(M\otimes \eta) < 2$, then 
$h^0(M\otimes \eta) = 0$. By Theorem~\ref{h0istwo}, the multiplicity of the singular
point of $\Xi$ is $2$ in this case as well.

Suppose finally that $h^0(M) =
3$. Then $h^0(M\otimes \eta) = 1$. If $M$ has a base point, then $C$ is hyperelliptic,
again by Clifford's theorem, and so $\Xi$ does not have an
isolated singularity. Otherwise, $M$ defines a morphism from $C$ to $\Pee^2$, which is
birational since $M$ has prime degree and is an embedding since the genus of a smooth
plane quintic is $6$.  Thus the only case that can arise is: $C$ is a smooth plane
quintic,
$h^0(M) = 3$ and  $M$ is the pullback of $\scrO_{\Pee^2}(1)$, and  $h^0(M\otimes
\eta) = 1$. This concludes the proof.
\end{proof}

\begin{remark}  Smith and Varley have suggested another proof of Theorem~\ref{smoothcase}, using
the result of Varley \cite[Proposition 7(i)]{Varley}  that, if $M$ is a theta characteristic on
$C$ and
$L=\pi^*M$ lies on the Prym variety, then 
$\operatorname{mult}_L\Xi \equiv \operatorname{mult}_M\Theta \mod 2$, where $\Theta$ is the
canonical theta divisor of $J(C)$. Thus, if $L$ corresponds to a triple point of $\Xi$,
$h^0(M)$ is odd, by the Riemann singularity theorem, and hence is either $1$ or $3$. The
argument then proceeds along the lines of the proof of Theorem~\ref{smoothcase}.
\end{remark}

\section{Singular curves}

\subsection{Preliminaries}

We extend the results of the previous section to the case of singular curves. Let
$\mathcal{R}_{g}$ denote the space of admissible covers constructed in \cite{Beau1},
such that the associated Prym morphism $P\colon \mathcal{R}_{g} \to \mathcal{A}_{g-1}$
is defined, dominant, and proper, and in particular it is surjective. A point of
$\mathcal{R}_{g}$ will be denoted by
$(\widetilde{C},
\pi)$, with corresponding involution $\tau$, and the corresponding Prym variety by
$P(\widetilde{C}, \pi)$. Since we are only concerned with principally polarized abelian
varieties whose theta divisors have an isolated singular point, and hence are
irreducible, we  will always assume that 
$(\widetilde{C}, \pi)$ satisfies condition $(*)$ of Beauville: the fixed points of
$\tau$ are exactly the singular points of $\widetilde{C}$, and at a fixed point the
local branches of $\widetilde{C}$ are not exchanged. (Otherwise by \cite[Theorem
5.4]{Beau1}, $P(\widetilde{C}, \pi)$ is either reducible or a Jacobian.) Moreover,  we
can also assume that, if there is a decomposition
$\widetilde{C} = \widetilde{C}_1\cup \widetilde{C}_2$, with
$\widetilde{C}_1\cap \widetilde{C}_2$ finite, then $\#(\widetilde{C}_1\cap
\widetilde{C}_2)$ is even and $\geq 4$ since if  $\#(\widetilde{C}_1\cap
\widetilde{C}_2)=2$, then,  by \cite[Lemma 4.11]{Beau1}, $P(\widetilde{C}, \pi)$ is
reducible. Of course, since the components and double points of $\widetilde{C}$ are
identified with those of $C$, a similar statement holds for $C$.

 We will need the
following numerology:

\begin{lemma}\label{numerology} Let $L$ be a line bundle on $\widetilde{C}$ such that
$\operatorname{Norm}L =
\omega_C$, and hence $L\otimes \tau^*L =\omega_{\widetilde{C}}$. Suppose that the
irreducible components of $\widetilde{C}$ are $C_i, \dots, C_k$, with $p_a(C_i) = p_i$.
 Suppose that $\overline{C}_i$ is
the component of $C$ lying under $C_i$, with $p_a(\overline{C}_i) =\overline{p}_i$. Let
$b_i = \#(C_i \cap (\bigcup_{j\neq i}C_j))$. Then
$$\deg (L|C_i) = p_i -1+\frac{b_i}{2} = 2\overline{p}_i -2 +b_i.$$
In particular $\deg (L|C_i)$ is even and, as $b_i \geq 4$,   $\deg (L|C_i)\geq
2\overline{p}_i + 2\geq 2$. Finally, $p_i\geq 1$, since $C_i\to \overline{C}_i$ is a
double cover with at least
$4$ branch points. 
\qed
\end{lemma}

\subsection{Clifford's theorem and hyperelliptic curves}

\begin{definition} Let $C$ be a nodal curve. A $g_2^1$ on $C$ is a line bundle $L$ such
that $\deg L = 2$ and $h^0(C; L) = 2$. We call $L$ \textsl{nonsingular} if the base
locus of $L$ does not contain a singular point of $C$, i.e.\ if for all $x\in C_{\rm
sing}$, there exists an $s\in H^0(L)$ such that $s(x)\neq 0$.  Finally, $C$ is \textsl{hyperelliptic} if $p_a(C) \geq
2$, $C$ is stable, and there is a
finite degree two morphism $C\to \Pee^1$.
\end{definition}

\begin{remark}\label{remark33} (i) Suppose that $C$ is stable and that $L$ is a
nonsingular
$g^1_2$ on
$C$. If $L$ has a base point $x$, necessarily a smooth point of $C$, then $L\otimes
\scrO_C(-x)$ has degree one and
$h^0=2$, and thus defines a morphism $C\to \Pee^1$ of degree one which is an
isomorphism on one component and maps all others to points. Thus, every positive
dimensional fiber defines a decomposition $C=C_1\cup C_2$ with $\#(C_1\cap
C_2) =1$. In particular, $C$ has a separating node $p$, i.e.\ a point $p\in  C_{\rm
sing}$ such that $C-\{x\}$ is disconnected.

\noindent (ii) Suppose that $C$ is stable and that $L$ is a base point free $g^1_2$ on
$C$. If the corresponding morphism $C\to \Pee^1$ is not finite, and $C_1$ is a
one-dimensional connected component of a fiber, then $C=C_1\cup C_2$ with $\#(C_1\cap
C_2) \leq 2$. 
\end{remark}

We will need the following version of Clifford's theorem for singular curves:

\begin{proposition}\label{Cliff} Let $C$  be a nodal curve and $L$ a line bundle on $C$
of degree
$d$. Suppose that, for each irreducible component $C_i$ of $C$, there exist  sections of $L$ and also of
$L^{-1}\otimes \omega_C$ which do not vanish identically on $C_i$, and hence $0\leq \deg
L|C_i \leq \deg \omega_C|C_i$  for every $i$. Then 
$$h^0(C; L) \leq \frac{d}{2} + 1.$$
Moreover, if $C$ does not have a separating node, then equality holds  only if
$L=\scrO_C$, $L=\omega_C$, or $C$ has a nonsingular $g^1_2$.
\end{proposition}
\begin{proof} Note that $L$ satisfies the hypotheses of the proposition if and only if 
$L^{-1}\otimes \omega_C$ does, and since $\deg(L^{-1}\otimes \omega_C) = 2p_a(C) -2 -
d$ and $h^0(L) = d- p_a(C) + 1 + h^0(L^{-1}\otimes \omega_C)$, 
$L$ satisfies the conclusions of the theorem if and only if $L^{-1}\otimes \omega_C$
does.

The proof is by induction on the number of nodes of $C$. The case where there are no
nodes is the usual form of Clifford's theorem for a special line bundle $L$. In the
general case, if both $L$ and $L^{-1}\otimes \omega_C$ are nonsingular, then the result
is proved by Beauville \cite[Lemma 4.7]{Beau1}. Suppose that either $L$ or
$L^{-1}\otimes \omega_C$ fails to be nonsingular. By the symmetry between $L$ and
$L^{-1}\otimes \omega_C$, we may suppose that $L$ is singular.

Let
$x\in C_{\rm sing}$ be such that every section of $H^0(L)$ vanishes at $x$. Let
$\nu\colon \hat C\to C$ be the result of normalizing $C$ at $x$, with $\nu^{-1}(x)
=\{p,q\}$. By assumption, $\nu^*H^0(L) \subseteq H^0(\hat C; \nu^*L(-p-q))$, and in
fact $\nu^*H^0(L) = H^0(\hat C; \nu^*L(-p-q))$, so that, for every component $D$
of $\hat C$, there exists a section of $\nu^*L(-p-q)$ which does not vanish
identically on
$D$. Set
$\hat L =
\nu^*L(-p-q)$.  Since $\omega_{\hat C} =
\nu^*\omega_C(-p-q)$, $(\hat L)^{-1}\otimes \omega_{\hat C} = \nu^*(L^{-1}\otimes
\omega_C)$. Thus a section of $L^{-1}\otimes \omega_C$ which does not vanish on any
component of $C$ defines a section of  $(\hat L)^{-1}\otimes \omega_{\hat C}$ with the
same property, so that $\hat L$ satisfies the inductive hypothesis. 

\noindent \textbf{Case 1}: $\hat C$ is connected, i.e.\ $x$ is not a separating node.
Then by the inductive hypothesis,
$$h^0(L) = h^0(\hat C;\hat L) \leq \frac{d-2}{2} + 1 = \frac{d}{2} < \frac{d}{2} + 1.$$
Hence strict inequality holds in this case.

\noindent \textbf{Case 2}: $\hat C$ is  not connected, i.e.\ $x$ is  a separating
node. Let $C_1, C_2$ be the connected components of $\hat C$. Let $L_i = \hat
L|C_i$ and let $d_i = \deg (L|C_i )$. Then $\deg L_i = d_i -1$ and $d=d_1+d_2$. Again,
$\hat L$ satisfies the hypotheses of Clifford's theorem, so that, by the inductive
hypothesis, 
$$h^0(L) = h^0(C_1; L_1) +h^0(C_2; L_2) \leq \frac{d_1-1}{2} + 1+\frac{d_2-1}{2} + 1 =
\frac{d}{2} + 1.$$
Thus the conclusion holds in this case as well. Moreover, if there is no separating
node, then the proof of Case I shows that equality holds only in the nonsingular case,
and then by Beauville's result either $L$ is $\scrO_C$ or $\omega_C$ or $C$ has a
nonsingular $g^1_2$.
\end{proof}  

\begin{lemma}\label{below} Let $\pi \colon \widetilde{C} \to C$ be an admissible double
cover satisfying condition $(*)$ of Beauville. If $\widetilde{C}$ has a nonsingular
$g^1_2$, then  $C$ has a nonsingular $g^1_2$.
\end{lemma}
\begin{proof} We shall just write out the proof in the case where the nonsingular
$g^1_2$ defines a finite morphism $f\colon \widetilde{C} \to \Pee^1$, the only case we
shall need. In this case, if $p$ and $q$ are smooth points of $\widetilde{C} $ such
that $f(p) = f(q)$, then in the usual way $h^0(\omega_{\widetilde{C}}(-p-q)) =
h^0(\omega_{\widetilde{C}}) -1$, i.e.\ a section of $\omega_{\widetilde{C}}$ which
vanishes at $p$ also vanishes at $q$. If $\pi(p)\neq \pi(q)$, since $\pi^*\omega_C =
\omega_{\widetilde{C}}$, it follows that every section of $\omega_C$ vanishing at
$\pi(p)$ vanishes at $\pi(q)$ and $C$ is hyperelliptic. Otherwise, for every choice of
smooth points $p$ and $q$ as above, if $f(p) = f(q)$, then $\pi(p) = \pi(q)$. But then
$f$ factors through $C$ and induces a finite, degree one morphism from $C$ to $\Pee^1$,
which is impossible.
\end{proof}

\subsection{Invariant line bundles}

We begin with the following situation: $\pi\colon \hat C \to N$ is a smooth
connected double cover, corresponding to the involution $\tau$ on $\hat C$. Suppose that
the fixed points of $\tau$ are $p_1, \dots, p_{2k}$ where $k>0$. Let $L$ be a line
bundle on
$\hat C$ such that $\tau^*L=L$ in $\Pic\hat C$, and let $(\Pic\hat C)^\tau$ be the
subgroup of all invariant line bundles. We define an invariant
$\alpha(L)$ as follows: choose an isomorphism $\varphi\colon \tau^*L \to L$ such that
$\varphi^2 =
\Id$, where $\varphi^2$ is understood to be the composition 
$$L = (\tau^*)^2L \xrightarrow{\tau^*\varphi} \tau^*L \xrightarrow{\varphi} L.$$
Such a $\varphi$ is unique up to $\pm 1$. For a fixed point $p_i$, $\varphi_{p_i}
\colon L_{p_i} \to L_{p_i}$ is multiplication by $\pm 1$. Thus
$$\alpha(L) = (\varphi_{p_1}, \dots, \varphi_{p_{2k}}) \in (\{\pm 1\})^{2k}/(\pm \Id)$$
is well-defined, where $\Id = (1, \dots, 1)$ and similarly for $-\Id$. Define the 
homomorphism $d\colon (\{\pm 1\})^{2k}
\to
\{\pm 1\}$   by
$d(\epsilon_1,\dots, \epsilon_{2k}) = \epsilon_1\cdots \epsilon_{2k}$. Note that $d$
factors through the quotient $(\{\pm 1\})^{2k}/(\pm \Id)$. The induced homomorphism
$(\{\pm 1\})^{2k}/(\pm \Id) \to \{\pm 1\}$ will again be denoted by $d$.

\begin{proposition} With notation as above, $\alpha \colon (\Pic \hat C)^{\tau} \to
(\{\pm 1\})^{2k}/(\pm \Id)$ is a surjective homomorphism satisfying
$$d(\alpha(L)) = (-1)^{\deg (L)}.$$
If $1\leq i_1<  \dots <  i_\ell\leq 2k$, then $\alpha(\scrO_{\hat C}(p_{i_1} + \dots +
p_{i_\ell}))$ is represented by $(\epsilon_1,\dots, \epsilon_{2k})$, where $\epsilon _i
= -1$ if and only if $i=i_j$ for some $j$.  Finally, $\Ker \alpha = \pi^*\Pic N$.
\end{proposition}
\begin{proof} It is easy to check that $\alpha$ is a homomorphism. The natural
linearization $\varphi$ on $\scrO_{\hat C}(p_{i_1} + \dots +
p_{i_\ell})$ induced by $f\mapsto \tau^*f$, where $f$ is a meromorphic function on a
Zariski open subset of $\hat C$ with at worst simple poles along the $p_{i_j}$, has
the property that it induces $-1$ on the fibers over the $p_{i_j}$, since there exists
an anti-invariant local uniformizing parameter over such points, and is $1$ over the
remaining fixed points. In particular, $\alpha$ is surjective and $d(\alpha(\scrO_{\hat
C}(p_{i_1} +
\dots + p_{i_\ell})) = (-1)^\ell$. 

The line bundle $L\in \Ker \alpha$ if and only if there is a linearization $\varphi$
which is the identity over the fixed points. By descent theory, this condition is
equivalent to $L=\pi^*M$ for some $M\in \Pic N$. Thus the set of line bundles 
$\{\scrO_{\hat C}(p_{i_1} + \dots + p_{i_\ell})\}$ contains a set of
representatives for $(\Pic \hat C)^{\tau}/\pi^*\Pic N$, and these representatives  are
unique modulo the relation
$\scrO_{\hat C}(p_1 + \dots + p_{2k}) =\pi^*\scrO_N(B)$ for the  divisor class $B$
which defines the double cover $\hat C$, i.e.\ $\scrO_N(2B) \cong \scrO_N(q_1 + \dots +
q_{2k})$ where $q_i = \pi(p_i)$. The formula $d(\alpha(L)) =
(-1)^{\deg (L)}$ then holds for all $L$, since it holds for $L=\pi^*M\in \pi^*\Pic N$
since $\deg \pi^*M$ is even, and it also holds for a set of coset representatives of
$(\Pic
\hat C)^{\tau}/\pi^*\Pic N$.
\end{proof}

We will also need to look at orders of vanishing of invariant or anti-invariant
sections. 

\begin{lemma}\label{order} Let $\varphi$ be a linearization on $L$ and let $p$ be a
fixed point of
$\tau$. 
\begin{enumerate}
\item[\rm (i)] If $\varphi_p = 1$, then every anti-invariant section $s$ of $L$ vanishes
at $p$, and every invariant section $s$ of $L$  which vanishes at $p$ vanishes to order
at least two at $p$.
\item[\rm (ii)] If $\varphi_p = -1$, then every invariant section $s$ of $L$
vanishes at $p$, and every anti-invariant section $s$ of $L$  which vanishes at $p$
vanishes to order at least two at $p$.
\end{enumerate}
\end{lemma}
\begin{proof} First suppose that $\varphi_p = 1$. Let $\sigma_0$ be a local generating
section of $L$ at $p$. By replacing $\sigma_0$ with $\sigma_0 +
\varphi(\tau^*\sigma_0)$ and shrinking further to a $\tau$-invariant open set
$U$ containing
$p$, we may assume that $\varphi(\tau^*\sigma) = \sigma$. If $s = f\sigma$ for some
holomorphic function $f$ on $U$, then
$s$ is invariant if and only if $\tau^*f= f$ and $s$ is anti-invariant if and only if
$\tau^*f= -f$. In the second case we see that $f$ vanishes at $p$, and in the first
case, if $f$ vanishes at $p$, then it vanishes to order at least $2$ there.

The case where $\varphi_p = -1$ follows from the case $\varphi_p = 1$ by replacing
$\varphi$ by $-\varphi$ .
\end{proof}

We now consider the case where $\pi\colon\widetilde{C}\to C$ is an admissible double
cover of the nodal curve $C$ satisfying $(*)$ and corresponding to the involution
$\tau$. Let
$\nu\colon
\hat C\to
\widetilde{C}$ be the (possibly disconnected) normalization of $\widetilde{C}$, so that
there is also the induced double cover $\hat C \to N$ and the corresponding involution,
which we continue to denote by $\pi$ and $\tau$. It is not always the case that a
$\tau$-invariant line bundle $\hat L$ on $\hat C$ lifts to a $\tau$-invariant line
bundle on $\widetilde{C}$. The condition is as follows:

\begin{proposition}\label{sameonboth} Let $\hat L\in (\Pic \hat C)^\tau$. Then $\hat L$
lifts to a
$\tau$-invariant line bundle on $\widetilde{C}$ if and only if there exists a
linearization $\hat\varphi$ of $\hat L$ such that, for every double point $x$ of
$\widetilde{C}$, if $\nu^{-1}(x) =\{p,q\}$, then
$\hat\varphi_p=\hat\varphi_q$. Hence, if $L$ is a  $\tau$-invariant line bundle  on
$\widetilde{C}$, then the pullback of $L$ to $\hat C$ is of the form
$\pi^*\hat{M}\otimes
\scrO_{\hat C}(\sum_{i=1}^k(p_i+q_i))$, where $x_1, \dots, x_k$ are distinct nodes of
$\widetilde{C}$ and $\nu^{-1}(x_i) = \{p_i, q_i\}$.
\end{proposition}
\begin{proof} If $\tau^*L\cong L$, then we can choose a linearization $\varphi$ on $L$
which then pulls back to a linearization $\hat\varphi$ of $\hat L$ satisfying
$\hat\varphi_p=\hat\varphi_q$ for every pair of points $p,q\in \hat C$ which are
identified under $\nu$.
Conversely,   such a linearization $\hat\varphi$ clearly induces an isomorphism from
$\tau^*L$ to $L$. 
\end{proof}

\begin{remark} It is easy to see that, if one lift of $\hat L$ is $\tau$-invariant,
then all lifts are $\tau$-invariant.
\end{remark}

Next, we need to know how to compute $\pi_*L$ for a $\tau$-invariant $L$.

\begin{proposition}\label{image}  Suppose that $\widetilde{C}_{\rm sing} =\{x_1, \dots,
x_n\}$ and set $\nu^{-1}(x_i) = \{p_i, q_i\}$. Let
$L$ be a $\tau$-invariant line bundle such that $\nu^*L\cong \pi^*\hat M_1(p_1+ q_1 +
\cdots + p_k + q_k)$. Thus we can also write $\nu^*L\cong \pi^*\hat M_2(p_{k+1}+ q_{k+1}
+ \cdots + p_n + q_n)$. Let $\nu_1\colon N_1\to C$ be the normalization of
$C$ at
$\pi(x_1),
\dots, \pi(x_k)$ and let $M_1$ be the line bundle on $N_1$ obtained from $\hat
M_1$ by using the identifications  $L_{p_{k+1}} \cong L_{q_{k+1}}$, 
\dots, $L_{p_n} \cong L_{q_n}$. Let $\nu_2\colon N_2\to C$ be
the normalization of
$C$ at
$\pi(x_{k+1}),
\dots, \pi(x_n)$ and let $M_2$ be the line bundle on $N_2$ obtained from $\hat
M_2$ by using the identifications  $L_{p_1} \cong L_{q_1}$, 
\dots, $L_{p_k} \cong L_{q_k}$. Then
$$\pi_*L \cong \nu_1{}_*M_1\oplus \nu_2{}_*M_2,$$
where, possibly after replacing $\varphi$ by $-\varphi$, $\nu_1{}_*M_1$ is the
$+1$-eigenspace for the action of $\varphi$ and $\nu_2{}_*M_2$ is the $-1$-eigenspace.
Finally, if
$\operatorname{Norm}L =
\omega_C$, then
$M_1$ is a theta characteristic on
$N_1$ and $M_2$ is a theta characteristic on $N_2$.
\end{proposition} 
\begin{proof}
This is essentially a local computation at each singular point. Let $R =
\Cee[[x,y]]/(xy)$ and let $S= \Cee[[u,v]]/(uv)$, and view $R$ as included in $S$ via
$x=u^2$, $y= v^2$. There is a natural involution $\tau$ on $S$, defined by $\tau(u) =
-u$ and $\tau(v) = -v$, and $R$ is the ring of $\tau$-invariants.

First consider $S$ as a module over $R$. Then $S = R \oplus( u\Cee[[x]] \oplus
v\Cee[[y]])$. In invariant terms, this says that
$$\pi_*\scrO_{\widetilde{C}} = \scrO_C \oplus \scrO_N(-B),$$
where $B$ is the line bundle on $N$ which is the square root of the branch divisor on
$\hat C \to N$. Moreover, the first factor is the $+1$-eigenspace and the second is the
$-1$-eigenspace. More generally, by the projection formula, 
$$\pi_*\pi^*M = M \oplus (M\otimes \scrO_N(-B)),$$
where in the usual linearization the first factor is the $+1$-eigenspace and the second
is the
$-1$-eigenspace.

Now suppose that $L$ is a line bundle on $C$ such that $\hat L \cong \scrO_{\hat
C}(p+q)$. Locally, there exists a $c\in \Cee^*$ such that $L$ corresponds to the
$S$-submodule $X$ of 
$\displaystyle \Cee[[u]]\cdot \frac{1}{u} \oplus \Cee[[v]]\cdot \frac{1}{v}$ consisting
of all
$(f_1, f_2)$ such that, if  the leading term of $f_1$ is $a/u$, then that of $f_2$ is
$ca/v$, i.e. 
$$X = \{ a\alpha + g_1+ g_2: g_1\in \Cee[[u]], g_2\in \Cee[[v]],
a\in
\Cee\}, $$
where $\alpha = 1/u + c/v$, so that $\tau(\alpha) =-\alpha$.
As an
$R$-module, 
$$X = (\Cee[[x]] \oplus \Cee[[y]]) \oplus (\Cee\cdot \alpha \oplus u\cdot \Cee[[x]]
\oplus v\Cee[[y]]) = X_+ \oplus X_-,$$
say, where $X_+$ is $\tau$-invariant and $X_-$ is $\tau$-anti-invariant. More
intrinsically, $X_+ \cong \widetilde{R}$, the normalization of $R$, and $X_- \cong R$.
These two cases together say that, if $L$ is a line bundle such
that $\nu^*L =\scrO_{\hat C}(p_1 + q_1 + \cdot + p_k+q_k)$ with the natural
linearization, then  the
$+1$-eigenspace of
$\pi_*L$ is isomorphic to $\nu_1{}_*\scrO_{N_1}$. Likewise, if $\nu^*L
=(\pi^*\hat M_1)\otimes \scrO_{\hat C}(p_1 + q_1 +
\cdot + p_k+q_k)$ with the natural linearization, then the
$+1$-eigenspace of
$\pi_*L$ is isomorphic to $\nu_1{}_*M_1$. By symmetry, a similar statement holds for
the $-1$-eigenspace.

To see the final statement, let $\hat L_1$ be the pullback of $L$ to the
normalization $N_1$ of $C$ at $x_1, \dots, x_k$, and similarly for $\hat L_2$.Then since
$\operatorname{Norm}L =
\omega_C$,
$$\operatorname{Norm}\hat L_1= \nu_1^*\omega_C = \omega_{N_1}(\pi(p_1) + \pi(q_1)
+\cdots + \pi(p_k) + \pi(q_k)).$$ But we also have $\hat L_1 = \pi^*M_1(p_1+q_1+ \cdots
+ p_k + q_k)$, so that
$$\operatorname{Norm}\hat L_1=   M_1^{\otimes 2}(\pi(p_1) + \pi(q_1)
+\cdots + \pi(p_k) + \pi(q_k)).$$
Comparing, we see that $ M_1^{\otimes 2} = \omega_{N_1}$. A similar statement works for
$M_2$.
\end{proof}

We will also use the following:

\begin{lemma}\label{311} Let $L$, $N_i$ and $M_i$ be as in the previous lemma. Suppose
that, for every irreducible component $D$ of $C$, 
$\deg(L|\widetilde{D})$ is even, where $\widetilde{D}$ is the component of
$\widetilde{C}$ corresponding to $D$. Let
$x_1,
\dots, x_\ell$ be the points lying on $D$ where $\nu_1\colon N_1 \to C$ is not an
isomorphism, with $\nu_1^{-1}(x_i) =\{p_i, q_i\}$. Suppose that $m$ of the points $p_i, q_i$
lie on the component of $N_1$ corresponding to $\nu_1^{-1}(D)$. Then:
\begin{enumerate}
\item[\rm (i)] $m = 2k$ is even, and $\deg (M_1|D) = \frac12\deg (L|\widetilde{D})
-k$.
\item[\rm (ii)] If $N_i = N'\cup N''$, where both $N'$ and $N''$ are a union of
components of $N_i$ and $N'\cap N''$ is finite, then $\#(N'\cap N'')$ is even.
\item[\rm (iii)] Every connected component of $N_i$ is either smooth rational or
semi\-stable with arithmetic genus at least one.
\end{enumerate}
\end{lemma}
\begin{proof} (i) follows from the definitions and the fact that $d(\alpha(L)) = 1$.
(ii) is clear if, say,
$N'$ is irreducible, and  follows in the general case by a straightforward induction on
the number of components of $N'$. (iii) is then clear.
\end{proof}

\subsection{Genus six}

Our goal is now to prove:

\begin{theorem}\label{mainsing} Let
$\pi \colon \widetilde{C}\to C$ be an admissible  double cover of a stable curve $C$
with $p_a(C) =6$. If the theta divisor $\Xi$ of the associated Prym variety
$P(\widetilde{C},
\pi)$ has a unique singular point which is a triple point, then $C$ is a
plane quintic and $\pi\colon \widetilde{C}\to C$ is the double cover associated to an
odd theta characteristic.
\end{theorem}

To prove Theorem~\ref{mainsing}, we may assume that $\pi$ satisfies condition $(*)$
and that, for every decomposition
$\widetilde{C} = \widetilde{C}_1\cup \widetilde{C}_2$, with
$\widetilde{C}_1\cap \widetilde{C}_2$ finite, then $\#(\widetilde{C}_1\cap
\widetilde{C}_2)\geq 4$.  Let $\Xi_{\rm sing} =\{x\}$, so that $\tau(x) = x$. Let
$L\in
\Pic \widetilde{C}$ be the corresponding line bundle, so that $\deg L =10$,  $h^0(L)
=2n$ is even, 
$\operatorname{Norm}L=\omega_C$ and 
$\tau^*L =L$. Thus $L= \omega_{\widetilde{C}}\otimes L^{-1}$. Via the isomorphism
$H^0(L) = H^0(M_1)
\oplus H^0(M_2)$ of Proposition~\ref{image}, we can write
$$\dim H^0(L) = 2n = \dim H^0(M_1) + \dim H^0(M_2) = n_1 + n_2,$$
say, where $n_i =\dim H^0(M_i)$. Possibly after relabeling, we may assume that $n_1\geq
n_2$. Then Theorem~\ref{mainsing} follows from:

\begin{theorem}\label{cases} With the above notation,
\begin{enumerate}
\item[\rm (i)] $2n = h^0(L) \leq 6$, and equality holds only if $C$ is hyperelliptic.
In this case $\dim \Xi_{\rm sing} \geq 1$.
\item[\rm (ii)] Suppose that $h^0(L) = 4$.
\begin{enumerate}
\item[\rm (a)] The case $(n_1, n_2) = (4,0)$ does not arise.
\item[\rm (b)] If $(n_1, n_2) =  (3,1)$, then either $C$ is a plane quintic and
$L=\pi^*M$, where $M$ is the unique line bundle of degree $5$ with $h^0(M) =3$ on $C$,
or
$\dim
\Xi_{\rm sing} \geq 1$.
\item[\rm (c)] If $(n_1, n_2) =  (2,2)$, then $L$ corresponds to a point of
multiplicity $2$ of $\Xi$.
\end{enumerate}
\item[\rm (iii)] Suppose that $h^0(L) = 2$.
\begin{enumerate}
\item[\rm (a)] If $(n_1, n_2) =  (2,0)$, then $L$ corresponds to a point of
multiplicity $2$ of $\Xi$.
\item[\rm (b)] If $(n_1, n_2) =  (1,1)$, then $L$ corresponds to a smooth point of
$\Xi$.
\end{enumerate}
\end{enumerate}
\end{theorem}

We now begin the proof of Theorem~\ref{cases}.

\begin{lemma}\label{connected} Let $s$ be a nonzero element of $H^0(L)$, and suppose
that
$\widetilde{C} = C_s'\cup C_s''$, where $C_s'$ is the union of the irreducible
components of
$\widetilde{C}$ where $s$ vanishes identically,  $C_s''$ is the union of the
irreducible components of
$\widetilde{C}$ where the restriction of $s$ is not zero, and hence
$C_s'\cap C_s''$ is finite. Then:
\begin{enumerate}
\item[\rm (i)] $C_s''$ is connected.
\item[\rm (ii)] Let $C_s'\cap C_s'' =\{x_1,
\dots, x_b\}$ and let $L'' = (L|C_s'')(-\sum _ix_i)$. Then $s$ defines  sections of
$L''$ and
$(L'')^{-1}\otimes \omega_{C_s''}$ which do not vanish on any component of $C_s''$.
\end{enumerate}
\end{lemma}
\begin{proof}
It suffices to consider the case when $C'_s$ is not empty. If $D$ is a
connected component of $C_s''$, then the intersection of $D$ with the remaining
components of $\widetilde{C}$ is the same as the intersection of $D$ with $C_s'$, and
thus $\#(D\cap C_s') \geq 4$. Hence $\deg(L|D) \geq 4$. By Lemma~\ref{numerology}, $L$
has degree at least two on every irreducible component of $C_s'$.
Thus, if $C_s''$ has at least two connected components, the only possibility is that
$C_s''$ has exactly two connected components $D_1$ and $D_2$, $L|D_i$ has degree $4$ for
$i=1,2$,  $C_s'$ is irreducible, $\#(C_s'\cap C_s'') =8$, and  $\deg(L|C_s')=2$. But by
Lemma~\ref{numerology},  $\deg(L|C_s')=2  = 2\overline{p} -2 + \#(C_s'\cap
C_s'')$, where $\overline{p}$ is the arithmetic genus of the image of $C'_s$
in $C$. Thus $\#(C_s'\cap C_s'') =4-2\overline{p} \leq 4$, contradicting  $\#(C_s'\cap
C_s'') =8$. Hence
$C_s''$ is connected, proving (i). 

Clearly the section $s$ does not vanish on any component of $C_s''$. Since 
\begin{align*}(L|C_s'')^{-1} \otimes \omega_{C_s''} &= (L|C_s'')^{-1} \otimes
(\omega_{\widetilde{C}}|C_s'')(-\sum _ix_i) \\&= (L^{-1}\otimes
\omega_{\widetilde{C}}|C_s'')(-\sum _ix_i)= (L|C_s'')(-\sum _ix_i)=L'',
\end{align*}
it follows that $(L'')^{-1} \otimes \omega_{C_s''} = L''(\sum_ix_i)$. The section $s$
then defines a section of $L''(\sum_ix_i) = (L'')^{-1} \otimes \omega_{C_s''}$ which
does not vanish on any component.
\end{proof}

\begin{corollary} If there exists a component of $\widetilde{C}$ on which every section
$s$ of $L$ vanishes, then $h^0(L) \leq 2$.
\end{corollary}
\begin{proof} As in the previous lemma, write $\widetilde{C} = C'\cup C''$, where
$C'$ is the union of all components of $\widetilde{C}$ where all sections of $L$
vanish, or equivalently where a generic section of $L$ vanishes, and let $C'\cap C''
=\{x_1,
\dots, x_b\}$. Then  $(L|C'')(-\sum _ix_i)$ has degree at most $8-4 =4$ and satisfies
the hypotheses of Proposition~\ref{Cliff}. Thus $h^0(L) = h^0((L|C'')(-\sum _ix_i))
\leq 3$, and so
$h^0(L) \leq 2$.
\end{proof}

\begin{corollary} We have $h^0(L) \leq 6$, and equality
holds only if $\widetilde{C}$ and hence $C$ are hyperelliptic.
\end{corollary}
\begin{proof} This is immediate from the previous corollary, Remark~\ref{remark33},
Proposition~\ref{Cliff} and Lemma~\ref{below}.
\end{proof}

\begin{proposition}\label{nonzero} Suppose that $s_+\in H^0(L)$ is a nonzero invariant
section and that
$s_-\in H^0(L)$ is a nonzero anti-invariant section. Then there exists a component $D$
of $\widetilde{C}$ such that neither $s_+$ nor $s_-$  restricts to $0$ on $D$.
\end{proposition}
\begin{proof} Suppose to the contrary that, for every component $D$
of $\widetilde{C}$ such that $s_+$ restricts non-trivially to $D$, $s_-$ restricts to
$0$ on $D$ and vice-versa.  Let $C_{s_+}''$ be the union of the components of
$\widetilde{C}$ where $s_+$ restricts non-trivially, and similarly for $C_{s_-}''$. By
assumption,
$C_{s_+}''$ and $C_{s_-}''$ meet in  a finite number (possibly zero) of points, and by
Lemma~\ref{connected}, $C_{s_+}''$ and $C_{s_-}''$ are connected. Let
$s$ be a generic linear combination of $s_+$ and $s_-$, and let $C_s''$ be the union of 
the components of $\widetilde{C}$ where $s$ has non-trivial restriction. Clearly $C_s''
= C_{s_+}'' \cup C_{s_-}''$. Again by Lemma~\ref{connected}, $C_s''$ is connected. It
follows that $C_{s_+}''\cap C_{s_-}'' \neq \emptyset$.

Let $C'_s$ be the union of the components where $s$ restricts to $0$. Clearly, $C'_s$
is the intersection of the corresponding subschemes $C_{s_+}'$ and $C_{s_-}'$, and $C =
C_s' \cup C_{s_+}''\cup C_{s_-}''$, where the pairwise intersections of $C_s' $,
$C_{s_+}''$, and $ C_{s_-}''$ are finite. First assume that $C_s' \neq \emptyset$. Then
as $\deg (L|C_{s_+}'') \geq 4$, $\deg (L|C_{s_-}'') \geq 4$, and $\deg (L|C_s') \geq 2$,
we must have  $\deg (L|C_{s_+}'') =\deg (L|C_{s_-}'') = 4$, and  both $C_{s_+}''$
and $C_{s_-}''$ meet the remaining components of $\widetilde{C}$ in $4$ points, at
least one of which lies on $C_{s_+}''\cap C_{s_-}''$. Let
$C_{s_+}''\cap\overline{(\widetilde{C} - C_{s_+}'')} = \{x_1, x_2, x_3, p\}$, say, where
$p$ corresponds to a point of $C_{s_+}''\cap C_{s_-}''$, and let
$C_{s_-}''\cap\overline{(\widetilde{C} - C_{s_-}'')} = \{y_1, y_2, y_3, q\}$, say,
where $q$ corresponds to the point $p$. Let $\varphi$ be a choice of linearization for
$L$ and let $\hat\varphi_\pm$ be the corresponding linearizations on $C_{s_\pm}''$. By
Proposition~\ref{sameonboth},
$(\hat\varphi_+)_p = (\hat\varphi_-)_q$. It then follows from Lemma~\ref{order} that
either $s_+$ vanishes to order at least two at $p$ or $s_-$ vanishes to order at least
two at $q$. In particular, either $\deg (L|C_{s_+}'') \geq 5$ or $\deg (L|C_{s_-}'')
\geq 5$. This is a contradiction, and hence there exists a component $D$ of
$\widetilde{C}$ where both $s_+$ and $s_-$ restrict non-trivially.

Finally assume that $C_s'=\emptyset$. Then $\widetilde{C} = C_{s_+}''\cup C_{s_-}''$,
and $\#(C_{s_+}''\cap C_{s_-}'') \geq 4$. At each point $p$ of $C_{s_+}''\cap
C_{s_-}''$, viewed as a point of $C_{s_+}''$, $s_+$ vanishes, and a similar statement
holds at the corresponding point $q$ of $C_{s_-}''$. Moreover, again by
Lemma~\ref{order}, either $s_+$ vanishes to order at least two at $p$ or $s_-$
vanishes to order at least two at $q$. It follows that $\deg (L|C_{s_+}'') +  \deg
(L|C_{s_-}'') \geq 4(3) = 12$, which contradicts $\deg L =10$. Thus, again there exists 
a component $D$ of
$\widetilde{C}$ where both $s_+$ and $s_-$ restrict non-trivially.
\end{proof}

\begin{proof}[Completion of the proof of Theorem~\ref{cases}] We have already
dealt with the case $h^0(L) \geq 6$. Next suppose that $h^0(L) = 2$ and $(n_1, n_2) =
(2,0)$. In this case, with notation as in Proposition~\ref{image}, given $q\in C$ a
smooth point such that
$\pi^{-1}(q) =\{p,\tau(p)\}$, $H^0(L\otimes \scrO_{\widetilde{C}}(2p+2\tau(p)))^+ =
H^0(M_1\otimes
\scrO_{N_1}(2q))$. Since $h^0(M_1) =2$, we may choose $q$ so that $h^1(M_1(2q)) =
h^0(M_1(-2q))=0$. The method of proof of Lemma~\ref{2.5} and Theorem~\ref{h0istwo} in
the case of a singular curve shows that
$L$ corresponds to a point of $\Xi$ of multiplicity two. If  $h^0(L) = 2$ and $(n_1,
n_2) = (1,1)$, then up to scalars there are unique nonzero sections $s_+ \in H^0(L)^+$,
$s_-
\in H^0(L)^-$ and there exists a component $D$ of $\widetilde{C}$ such that neither
$s_+$ nor $s_-$ restricts to $0$ on $D$. After relabeling, we may assume that
$H^0(L)^+ = H^0(M_1)$ and that $H^0(L)^-= H^0(M_2)$. For a general point
$p\in D$, lying over $q\in C$, we can assume that $H^0(M_1(-q)) =H^0(M_2(-q)) = 0$.
Thus 
$$H^0(L(-p-\tau(p))) = H^0(M_1(-q)) \oplus H^0(M_2(-q)) = 0,$$
so that $p$ and $\tau(p)$ impose independent conditions on $H^0(L)$. It then follows
from Theorem~\ref{SV} that the multiplicity of $\Xi$ at the point corresponding to $L$
is one. (Alternatively, one could use Mumford's tangent cone condition as proved by
Beauville~\cite[(4.1)]{Beau1} in the case of an admissible double cover.)

Now suppose that $h^0(L) = 4$. If $(n_1, n_2) = (2,2)$, we claim that there exist points
$q_1, q_2\in C$, with preimages $\pi^{-1}(q_i) =\{p_i,
\tau(p_i)\}$, such that the points $p_1, \tau(p_1),p_2, \tau(p_2)$ impose independent
conditions on $H^0(L)$. Then, again by Theorem~\ref{SV}, it follows that the
multiplicity of $\Xi$  at the point corresponding to $L$
is  two. To find the points $q_1, q_2$, begin with a component $D_1$ of $\widetilde{C}$
such that there exist sections both in $H^0(L)^+$ and $H^0(L)^-$ whose restrictions to
$D_1$ are not zero, and choose $q_1$ to be a general point of the image of $D_1$ in $C$.
The   arguments above for the case $(n_1, n_2) = (1,1)$ show that 
$$H^0(L(-p_1-\tau(p_1))) = H^0(M_1(-q_1)) \oplus H^0(M_2(-q_1)),$$
where now $\dim H^0(M_i(-q_1)) =1$. Choose a nonzero section 
$$s_+ \in
H^0(L(-p_1-\tau(p_1)))^+ = H^0(M_1(-q_1)),$$ and choose $s_-$ similarly. In particular,
we can view $s_\pm$ as nonzero sections of $L$.  By Proposition~\ref{nonzero}, there
exists a component
$D_2$ (possibly equal to
$D_1$) of
$\widetilde{C}$ such that both $s_+$ and $s_-$ have nonzero restriction to $D_2$. Let
$q_2$ be a general point of the image of $D_2$ in $C$ and let $\pi^{-1}(q_2) =\{p_2,
\tau(p_2)\}$. Then  
$$H^0(L(-p_1-\tau(p_1)-p_2-\tau(p_2))) = H^0(M_1(-q_1-q_2)) \oplus
H^0(M_2(-q_1-q_2))=0.$$
Hence $p_1, \tau(p_1),p_2, \tau(p_2)$ impose independent
conditions on $H^0(L)$, so that $L$ corresponds to a point of multiplicity two.

Thus we may henceforth assume that $(n_1, n_2) = (3,1)$ or $(4,0)$. It will be
convenient to divide the rest of the argument into cases.

\smallskip
\noindent\textbf{Case I}:   \emph{$N_1$ is connected, and, for every irreducible
component 
$D$ of $N_1$, not all sections of
$M_1$ vanish on $D$}. 

In this case, $M_1^{\otimes 2} =
\omega_{N_1}$, 
$M_1$ has degree at most
$5$, and, by Lemma~\ref{311}, $\deg M_1=5$  if and only if $N_1 = C$. By
Proposition~\ref{Cliff}, $h^0(M_1) \leq 3$. Moreover, if equality holds, either
$\deg M_1 = 5$ or $\deg M_1 =4$ and $C$ is obtained
from $N_1$ by identifying two points $p$ and $q$, which must then lie in the same 
component of $N_1$ by parity. Suppose that $\deg N_1=4$.  Then $N_1$ does not
have a separating node and thus, since equality holds in Clifford's theorem, 
$N_1$ has a nonsingular
$g^1_2$.  Since every irreducible component of $N_1$ meets the remaining components in
at least
$4$ points,
$N_1$ is stable and  hyperelliptic. Thus
$C$ satisfies (b) of \cite[Theorem 4.10]{Beau1}. In this case, the Prym variety is a
Jacobian, by the main theorem in \cite{Shok} (where these curves are called
\textsl{quasi-trigonal}), and hence $\dim \Xi_{\rm sing} \geq 1$. 

We may thus assume  that $\deg M_1 =5$, and hence that $N_1 = C$. If $M_1$ has a
base point which is a node, necessarily non-separating, then the proof of
Proposition~\ref{Cliff} shows that
$h^0(M_1)
\leq 2$, so that this case does not arise. If
$M_1$ has a base point which is a smooth point of
$C$, but no singular base point, then by Clifford's theorem $N_1 =C$ has a nonsingular
$g^1_2$ and so is hyperelliptic. Again,
$\dim \Xi_{\rm sing} \geq 1$. Otherwise  $M_1$ defines a finite morphism $f\colon C\to
\Pee^2$. If $f$ has positive degree on a component, or maps two components to
the same irreducible curve in $\Pee^2$, then there exists a point $x\in f(C)$ such that 
$f^{-1}(x) =\{y_1, \dots, y_k\}$ with $k\geq 2$, $x$ is a smooth point of
$f(C)$,  and each
$y_i$ is a smooth point of
$C$. It follows that $M_1(-\sum_iy_i)$ has degree at most three, is base point
free,  and
$h^0(M_1(-\sum_iy_i))=2$. The projection from $x$ defines a morphism $C\to \Pee^1$,
which is finite by the arguments of Remark~\ref{remark33}. Thus
$C$ is either trigonal or hyperelliptic. By the main result of \cite{Shok}, $\dim
\Xi_{\rm sing} \geq 1$ in this case as well. In the remaining case, the morphism $f$
satisfies:
$f_*\scrO_{C}/\scrO_{f(C)}$ is a skyscraper sheaf. Since $f(C)$ is a plane quintic,
$p_a(C) = p_a(f(C)) = 6$, and hence
$f_*\scrO_{C}=\scrO_{f(C)}$. It follows that $f$ is an embedding, 
$C$ is a plane quintic and $M_1$ is the pullback of $\scrO_{\Pee^2}(1)$.

\smallskip
\noindent\textbf{Case II}:   \emph{$N_1$ is connected, but there is a component of
$N_1$ on which all sections of
$M_1$ vanish}. 

By Lemma~\ref{311}, either $N_1$ is semistable or it is smooth rational,
in which case $h^0(N_1) =0$. Assume that $N_1$ is semistable and let $\overline{N_1}$ be
the associated stable curve, so that the theta characteristic $M_1$ is the pullback of
a theta characteristic $\overline{M_1}$ on $\overline{N_1}$ with $h^0(\overline{M_1})
= h^0(M_1)$. Note that  $\overline{M_1}$ has strictly positive degree on every
component of $\overline{N_1}$. Clearly
$\overline{N_1}$ and $\overline{M_1}$ satisfy the conclusions of Lemma~\ref{311}. If 
the only components of $N_1$ where all sections  of $M_1$ vanish are
contained in  fibers of the morphism
$N_1\to \overline{N_1}$, then, for every irreducible component $\overline{D}$ of 
$\overline{N_1}$, not all sections of 
$\overline{M_1}$ vanish on $\overline{D}$, but there is a node of $\overline{N_1}$, the
image of a component where all sections of $M_1$ vanish, which is in the
base locus of $\overline{M_1}$. It follows  as in Case I that
$h^0(\overline{M_1})\leq 2$. Thus, we may assume that there is a component of
$\overline{N_1}$ where all sections of
$\overline{M_1}$ vanish. Write $\overline{N_1} = \overline{N}'\cup \overline{N}''$,
where
$\overline{N}'\neq
\emptyset$ is the union of all components of
$\overline{N}_1$ on which all sections of $\overline{M_1}$ vanish. Note that 
$\#(\overline{N}'\cap \overline{N}'')$
is even and $\deg(\overline{M_1}|\overline{N}'')= \deg \overline{M_1} -
\deg(\overline{M_1}|\overline{N}') \leq 5-1=4$. If
$\overline{N}'\cap
\overline{N}''=\{x_1, \dots, x_\ell\}$,  where the $x_i$ are distinct and $\ell
=\#(\overline{N}'\cap \overline{N}'')$, then
$h^0(\overline{M_1}) = h^0((\overline{M_1}|\overline{N}'')(-\sum _ix_i))$, where 
$$0\leq \deg(\overline{M_1}|\overline{N}'')(-\sum _ix_i) = \deg
(\overline{M_1}|\overline{N}'')-\ell.$$
Thus $\ell \leq 4$, so that $\ell$ is either $2$ or $4$. If
$\ell =\#(\overline{N}'\cap \overline{N}'')= 4$, then 
$\deg(\overline{M_1}|\overline{N}'')(-\sum _ix_i)\leq 0$. Since each connected
component of $\overline{N}''$ meets the rest of $\overline{N_1}$ in at least two
points,  there are at most two connected components of
$\overline{N}''$ and so
$h^0(\overline{M_1}) = h^0((\overline{M_1}|\overline{N}'')(-\sum _ix_i)) \leq 2$. If
$\ell =2$, then $\overline{N}''$ is connected and the morphism
$\overline{N}''\to C$ must have involved a normalization of $C$ at at least two nodes.
In this case,
$\deg (\overline{M_1}|\overline{N}'')
\leq 3$, $\deg(\overline{M_1}|\overline{N}'')(-\sum _ix_i) \leq 1$, and hence
$$h^0(\overline{M_1}) = h^0((\overline{M_1}|\overline{N}'')(-\sum _ix_i)) \leq 2.$$
Hence neither case arises if $n_1$ is
$3$ or $4$.

\smallskip
\noindent\textbf{Case III}: \emph{$N_1$ is not connected}. 

Let $D_1, \dots, D_t$ be the
connected components of $N_1$ and let $C_i$ be the subscheme of $\widetilde{C}$
corresponding to $D_i$. If $d_i=\frac12\deg (L|C_i)$, then $d_i>0$ and $\sum _id_i =5$.
By Lemma~\ref{311}, $D_i$ is semistable or
$D_i$ is smooth rational, and $\deg (M_1|D_i) \leq d_i - 2$ since the image of $D_i$ in
$C$ must meet the remaining components in at least $4$ points. Moreover,
$M_1$ is a theta characteristic on every component $D_i$. Finally, each
irreducible component of $D_i$ meets the remaining components of $D_i$ in an even
number of points. The possibilities for
$(d_1,
\dots, d_t)$ are, after reordering, $(4,1)$, $(3,2)$, $(3,1,1)$, $(2,2,1)$, $(2,1,1,1)$,
$(1,1,1,1,1)$, and thus $\deg(M_1|D_i) \leq 2$. If
$\deg(M_1|D_i) \leq -1$, then $D_i$ is
smooth rational and $h^0(M_1|D_i) = 0$. If $\deg(M_1|D_i) =0$, then $D_i$ is
semistable of arithmetic genus one, and
$h^0(M_1|D_i) \leq 1$. If $\deg(M_1|D_i) =1$, then $d_i\geq 3$, $D_i$ is
semistable of arithmetic genus two, there is exactly one irreducible component of
$D_i$ such that $M_1$ restricts to a line bundle of degree one, and the remaining
components are smooth rational curves meeting the rest of $D_i$ in two points. In
particular, $M_1|D_i$ is induced by a line bundle $\overline{M_1}$ of degree one on
the stable model
$\overline{D_i}$ of $D_i$, which is irreducible of arithmetic genus two. It follows 
that $h^0(M_1|D_i)=h^0(\overline{M_1}) \leq 1$. If 
$D_i$ is a component such that $\deg(M_1|D_i)=2$, then $d_i =4$ and $M_1|D_i$ is induced
by a line  bundle $\overline{M_1}$ of degree two on the stable model
$\overline{D_i}$ of $D_i$, which has arithmetic genus three and at most two components.
If there exist sections of $\overline{M_1}$ which do not vanish on a component,
then  Proposition~\ref{Cliff} implies that  $h^0(M_1|D_i)=h^0(\overline{M_1}) \leq 2$.
Otherwise, it is easy to see that  $h^0(M_1|D_i) =h^0(\overline{M_1})=0$. Running
through the list of possible $(d_1,
\dots, d_t)$ above, it is easy to see that in all cases $h^0(M_1) = \sum
_{i=1}^th^0(M_1|D_i) \leq 2$.

 Thus, examining all of the cases, we see that $h^0(M_1) =4$ is impossible, and
if $h^0(M_1) =3$, then either $N_1=C$ and $M_1$ embeds $C$ as a plane quintic, or  $C$
is trigonal or hyperelliptic or is obtained from a hyperelliptic curve by identifying two
points. Thus, if $C$ is not a plane quintic, then 
$\dim
\Xi_{\rm sing} \geq 1$. This concludes the proof of Theorem~\ref{cases}.
\end{proof}

\section{The proof of Theorem~\ref{ourthm}}

\subsection{Plane quintics and cubic threefolds}

Let $C$ be a reduced nodal  curve, with normalization $\nu\colon N \to C$, and let
$M$ be a theta characteristic on $N$. If the preimages of the nodes of $C$ are $p_1, q_1,
\dots, p_\delta, q_\delta$, then a local calculation shows that
\begin{align*}
Hom(\nu_*M, \omega_C) &=\nu_*\left(M^{-1}\otimes \scrO_N\left(-\sum_i(p_i+q_i)\right)
\otimes
\nu^*\omega_C\right)\\& =
\nu_*(M^{-1}\otimes K_N) = \nu_*M.
\end{align*}
Using this, we have the following result of Beauville \cite[Prop. 4.2]{Beau3} (cf.\
also \cite[6.27]{Beau2} as well as \cite[Proposition V.2.3]{DS}):

\begin{theorem}\label{Bthm} Let $C$ be a nodal plane quintic curve and $M$ a theta
characteristic on the normalization $N$ of $C$ such that $h^0(M) = 1$. Then there exists a
symmetric matrix 
$$A =\begin{pmatrix} \ell_1 &\ell_2 &q_1\\
\ell_2 &\ell_3 & q_2\\
q_1 & q_2 & f\end{pmatrix},$$
where the $\ell_i$ are linear polynomials in $x,y,z$, the $q_i$ are quadratic, and $f$ is
cubic, such that the following is exact:
$$0 \to \scrO_{\Pee^2}(-2)^2 \oplus \scrO_{\Pee^2}(-3) \xrightarrow{A} \scrO_{\Pee^2}(-1)^2
\oplus \scrO_{\Pee^2} \to \nu_*M \to 0.$$
In particular $\det A =\Delta$ is a polynomial of degree $5$ defining $C$. \qed
\end{theorem}

With hypotheses as in Theorem~\ref{Bthm}, let $x,y,z,w,t$ be homogeneous coordinates on
$\Pee^4$ and consider the cubic threefold $X$ defined by  the homogeneous cubic polynomial
$$F(x,y,z,w,t) = \ell_1w^2 + 2\ell_2wt + \ell_3t^2 + 2q_1w+ 2q_2t + f.$$

\begin{proposition}\label{cubicexists} The cubic threefold $X$ contains the line
$\lambda$ defined by
$x=y=z=0$. Under the assumption that $C$ is nodal and that $M$ is a theta
characteristic on $N$ with
$h^0(M) =1$, the threefold $X$ is smooth. Finally, the curve $C$ together with the theta
characteristic
$M$ on the normalization $N$ is the discriminant curve of the pair $(X, \lambda)$.
\end{proposition}
\begin{proof} Clearly $X$ contains $\lambda$. To see that $X$ is smooth, assume that $p_0$
is a point of $\operatorname{Sing}X$. We first consider the case
$p_0=(x_0, y_0, z_0, w_0, t_0)$ with $(x_0, y_0, z_0) \neq (0,0,0)$,
i.e.\ $p_0\notin \lambda$. After a change of coordinates fixing $x,y,z$, we may assume that
$(w_0, t_0) = (0,0)$. Since $p_0$ is a singular point of $X$,  
$$f(x_0, y_0, z_0) = q_1(x_0, y_0, z_0) = q_2(x_0, y_0, z_0) =0$$ 
and $(x_0, y_0, z_0)$ is a
singular point of $f$. In particular, the matrix
$A$ is not invertible at
$(x_0, y_0, z_0)$.  Let $R$ be the local ring of $\Pee^2$ at $(x_0, y_0, z_0)$, with 
maximal ideal $\mathfrak{m}$, and continue to denote by $A$ the matrix with coefficients in
$R$ induced by
$A$.  First assume 
that the matrix $A\mod \mathfrak{m}$ has rank two at
$(x_0, y_0, z_0)$. Then it is easy to check that
$R^3/A\cdot R^3\cong R/\Delta$. After a linear change of coordinates in $w,t$, we may
assume that
$A\mod
\mathfrak{m}$ is  of the form
$\displaystyle\begin{pmatrix} 1&0&0\\0&1&0 \\0&0&0 \end{pmatrix}$.
Thus, $\Delta \equiv f\mod \mathfrak{m}^2$.  Since $(x_0, y_0, z_0)$ is a singular point of
the cubic defined by $f=0$, $\Delta$ has a singular point at $(x_0, y_0,
z_0)$, and the cokernel of $A$ is locally of the form $\scrO_C$, contradicting the
hypotheses on $A$.

Next suppose that $A\mod \mathfrak{m}$ has rank one. Then we may assume that 
$A\mod \mathfrak{m}$ is  of the form
$\displaystyle\begin{pmatrix} 1&0&0\\0&0&0 \\0&0&0 \end{pmatrix}$.
Note that
$\Delta \equiv \ell_3f- q_2^2 \mod \mathfrak{m}^3$. Since $f$ is singular at $(x_0, y_0, z_0)$, 
$\ell_3f- q_2^2 \equiv -q_2^2
\mod \mathfrak{m}^3$, and hence the tangent cone of $\ell_3f- q_2^2$ is degenerate. This 
contradicts the hypothesis on $\Delta$.

The remaining case is where the point $p_0$ lies on $\lambda$, in which case we may 
assume that
$p_0 = (0,0,0,1,0)$. Then $p_0$ is a singular point of $X$ if and only if $\ell_1=0$. In this
case,
$\Delta = -(q_1^2\ell_3 - 2\ell_2q_1q_2 + \ell_2^2f)$ automatically has a singular point at
$(x_0, y_0, z_0)$ lying in the intersection of $\ell_2 = q_1=0$. At such a point, if
$\ell_3f-q_2^2\neq 0$, then the matrix $\displaystyle \begin{pmatrix}
\ell_3&q_2\\q_2&f\end{pmatrix}$ is nonsingular, and it is easy to see that  $R^3/A\cdot
R^3
\cong  R/\Delta$, a contradiction. Otherwise, it is easy to check that $\Delta$ is 
congruent to a square mod $\mathfrak{m}^3$, and so cannot have a nondegenerate tangent cone.
For example, if there exists a $c$ such that $\ell_3\equiv cq_2\mod \mathfrak{m}$ and
$q_2\equiv cf\mod \mathfrak{m}$, then 
$$\Delta \equiv -f(cq_1-\ell_2)^2 \mod \mathfrak{m}^3.$$
Once again, we reach a contradiction.

The final point to check is that the theta characteristic $M$ arises from the conic 
 bundle structure on $X$. Since $X$ is nonsingular, $h=\ell_1\ell_3-\ell_2^2$ is
nonzero. First assume that both $C$ and the conic $Q$ defined by $h$ are smooth. We
claim that $Q$  meets
$C$ everywhere with even multiplicity, and the support of the intersection is exactly
the support of the sheaf
$\nu_*M/\scrO_C$, where the inclusion $\scrO_C\to \nu_*M$ is defined by the unique
section of $\nu_*M$.  Note that 
$$\det A = fh + q_1(\ell_2q_2-q_1\ell_3) -q_2(\ell_1q_2-\ell_2q_1) = fh +(\ell_3q_1^2
-2\ell_1q_1q_2 +\ell_1q_2^2),$$
and hence, if $p\in \Pee^2$ and $R=\scrO_{\Pee^2,p}$, $$R/(\det A, h) = R/(h,
\ell_3q_1^2 -2\ell_1q_1q_2 +\ell_1q_2^2).$$
Since  $Q$ is smooth,  at least one of $\ell_1, \ell_3$ is nonzero
at $p$. If for example
$\ell_3\neq 0$, then 
$$\ell_3q_1^2 -2\ell_1q_1q_2 +\ell_1q_2^2 = \frac{1}{\ell_3}(\ell_3q_1-\ell_2q_2)^2.$$
The defining exact sequence for $\nu_*M$ shows
that, at $p$,  
$$(\nu_*M/\scrO_C)_p \cong R^2/((\ell_1, \ell_2), (\ell_2, \ell_3), (q_1,
q_2)).$$  The homomorphism $R^2\to R/h$ defined by $(a,b) \mapsto
a\ell_3-b\ell_2$ then defines an isomorphism from $\nu_*M/\scrO_C$ to $R/(h,
\ell_3q_1-\ell_2q_2)$. Thus  $\nu_*M/\scrO_C$ has
support at $p$ if and only if $p\in C\cap Q$, and the length of the support there is
one-half the local intersection number. 
 It follows that the function field of $\widetilde{C}$ is generated by
$\sqrt{\ell_1\ell_3-\ell_2^2}$. By
\cite[Lemma 1.6]{Beau2}, the same holds for the double cover of $C$ defined by the 
discriminant curve construction. Hence the two double covers agree. This completes the
proof when $C$ and $Q$ are smooth. For the general case, choose a pencil of cubic
threefolds $X_t$ containing the line $\lambda$, specializing to the given pair
$(X,\lambda)$, and such that,  for the general member of the pencil, the discriminant
curve and the conic defined by the first
$2\times 2$ minor of the matrix $A_t$ are smooth. A limiting argument,  using the fact
that $\mathcal{R}_6$ is separated, then imples that the corresponding admissible double
covers agree for all values of $t$.
\end{proof}

\subsection{Completion of the proof}

Let $A$ be a principally polarized abelian variety of dimension $5$. Then $A\cong
P=P(\widetilde{C},
\pi)$ as principally polarized abelian varieties, where $\pi\colon\widetilde{C}\to C$
is an admissible double cover. If the theta divisor $\Xi$ of $P$ has a unique singular
point which is a triple point, then, by Theorem~\ref{smoothcase} and
Theorem~\ref{mainsing}, $C$ is a plane quintic and the theta characteristic $M$ on the
normalization of $C$ corresponding to the double cover
$\pi\colon \widetilde{C}\to C$ is odd. By Proposition~\ref{cubicexists}, there exists a
smooth cubic threefold $X$  and a line $\lambda$ on $X$ such that the discriminant
curve of the pair $(X,\lambda)$ is $C$ and the corresponding  theta characteristic on
$N$ is 
$M$. By Beauville's result \cite[Proposition 2.8]{Beau2}, $JX \cong P$ as principally
polarized abelian varieties, and hence $A\cong JX$. This concludes the proof of
Theorem~\ref{ourthm}.

\bigskip
\noindent
Department of Mathematics \\
Columbia University \\
New York, NY 10027 \\
USA

\bigskip
\noindent
{\tt  casa@math.columbia.edu, rf@math.columbia.edu}

\end{document}